\documentclass[12pt,leqno]{amsart}
\usepackage{euscript}
\newtheorem{Appendix}{A.\hskip -.01em}
\newtheorem{theorem}{Theorem}
\newtheorem{definition}[theorem]{Definition}
\newtheorem{corollary}[theorem]{Corollary}

\newtheorem{example}[theorem]{Example}

\newtheorem{proposition}[theorem]{Proposition}

\newtheorem{problem}[theorem]{Problem}
\newtheorem{conjecture}[theorem]{Conjecture}

\newcommand{\Dtriangle}[6]{
\setlength{\unitlength}{.6cm}
\begin{picture}(6,4.1)(0,-.1)
\thinlines
\put(0,3){\makebox(0,0){$#1$}}
\put(6,3){\makebox(0,0){$#2$}}
\put(0,0){\makebox(0,0){$#3$}}

\put(1,3){\vector(1,0){4}}
\put(0,2.5){\vector(0,-1){2}}
\put(1,0.5){\vector(2,1){4}}

\put(-.5,1.5){\makebox(0,0)[r]{\scriptsize $#5$}}
\put(3,3.5){\makebox(0,0)[b]{\scriptsize $#4$}}
\put(3.5,1.25){\makebox(0,0)[l]{\scriptsize $#6$}}
\end{picture}
}

\def\Cact{{\mathcal C} \hskip -.1em{\it act\/}} \def\disc{{\EuScript D}}
\def\calS{{\mathcal S}}  \def\calD{{\mathcal D}} 
\def\ccc{{\it cup\/}} \def\Braid{{\mathcal B} \hskip -.1em{\it raid\/}}

\def\susp{\hskip .2em \uparrow \hskip -.2em} \def\calA{{\mathcal A}}
\def\ut{{\underline t}} \def\desusp{\downarrow} \def\Span{{\it Span\/}}
\def\Gerst{{\mathcal G}\hskip -.1em {\it er\/}} 
\def\BbbZ{{\mathbb Z}} \def\Reg{{\it Reg\/}}
\def\calB{{\mathcal B}} \def\calP{{\EuScript P}}  
\def\Endom{{\it End\/}} \def\calZ{{\mathcal Z}}
\def\CPAA#1{{\mathcal C}_{\EuScript P}^{#1}(A;A)} 
\def\bfk{{\bf k}}\def\id{1\!\!1} \def\calQ{{\EuScript Q}}
\def\Ass{{\mathcal A{\it ss\/}}} \def\Com{{\mathcal C{\it om\/}}}
\def\Lie{{\mathcal L{\it ie\/}}} \def\uAss{{\underline{\mathcal A{\it ss\/}}}}
\def\preLie{{{\it pre}\mathcal L{\it ie\/}}} \def\Ker{{\it Ker\/}} 
\def\Sym{{\mathcal S{\hskip -.1em \it ym\/}}} 
\def\Mag{{{\EuScript M}\hskip -.05em {\it ag\/}}}
\def\uMag{{\underline{{\EuScript M}\hskip -.05em {\it a}}{\it g\/}}}
\def\Lin{{\it Lin\/}} 
\def\Lie{{{\mathcal L}{\it ie}}} \def\ot{\otimes} \def\sgn{{\it sgn\/}}

\def\qed{\hspace*{\fill}\mbox{\hphantom{mm}\rule{0.25cm}{0.25cm}}\\}
\def\cases#1#2#3#4{
                  \left\{
                         \begin{array}{ll}
                           #1,\ &\mbox{#2}
                           \\
                           #3,\ &\mbox{#4}
                          \end{array}
                   \right.
}
\def\HAP#1{H^{#1}_{\calP}(A;A)}\def\CAD#1{C^{#1}_{\calD}(A;A)}
\def\CAP#1{C^{#1}_{\calP}(A;A)}  \def\As{\Ass}
\def\otexp#1#2{#1^{\otimes #2}} \def\tildeCAP#1{\widetilde C^{#1}_{\calP}(A;A)}
\def\End{{\mathcal E} \hskip -.1em{\it nd\/}} \def\calT{{\EuScript T}}
\def\freepreLie{{\it pre\/}{\mathbb L}} 
\def\freeLie{{\mathbb L}} \def\ucalP{{\underline{\EuScript P}}}

\def\uchi{{\underline{\chi}}} \def\ns{non-$\Sigma$} 
\def\epi{\to \hskip -.6em \to} \def\ss{{\mathbf s \hskip .1em}} 
\def\opt{\calP \ot \calP^!} \def\triv{{\bf 1}} \def\Im{{\it Im\/}}

\begin{document}
\bibliographystyle{plain}
\baselineskip 1.3em

\title{Cohomology operations and the Deligne conjecture}
\author[M. Markl]{Martin MARKL}
\thanks{The author was supported by the grant GA \v CR 201/05/2117 and by
   the Academy of Sciences of the Czech Republic, 
   Institutional Research Plan No.~AV0Z10190503.}
\date{June~10, 2005}
\email{markl@math.cas.cz}

\begin{abstract}
The aim of this note, which raises more questions than it answers,  
is to study natural operations acting on the cohomology of various types of
algebras. It contains a lot of very surprising partial results and examples. 
\end{abstract}

\address{Mathematical Institute of the Academy, {\v Z}itn{\'a} 25, 
         115 67 Prague 1, The Czech Republic}

\maketitle

\section*{Introduction}

In this note, 
all algebraic objects will be defined over a fixed field $\bfk$ of
characteristic zero. 
An~{\em algebra\/} means an algebra over a quadratic Koszul
operad~$\calP$~\cite[II.3.3]{markl-shnider-stasheff:book}. 
This generality covers all ``reasonable'' algebras -- associative, Lie, 
commutative associative, Poisson, Gerstenhaber, Leibniz,~\&c.

By the {\em cohomology\/} of a $\calP$-algebra $A$ we mean the
operadic cohomology $H^*_\calP(A;A)$ of $A$ with coefficients in
itself~\cite[II.3.100]{markl-shnider-stasheff:book}, defined as the
cohomology of the cochain complex $\CAP* = (\CAP*,d_\calP)$ recalled
in A.\ref{Ferda} of the appendix to this note. The complex $\CAP*$
generalizes the ``standard constructions'' and $H^*_\calP(A;A)$ the
``classical'' cohomology (Hochschild for associative algebras,
Chevalley-Eilenberg for Lie algebras, Harrison for associative
commutative algebras,~\&c.)  In general, $H^*_\calP(A;A)$ agrees with
the triple cohomology~\cite[Proposition~8.6]{fox-markl:ContM97} and governs
deformations of $A$ in the category of $\calP$-algebras.

By a {\em natural operation\/} we mean a multilinear operation on
$\HAP *$ induced by a natural multilinear cochain operation on
$\CAP*$. Naturality means being defined using data that do not
depend on a concrete algebra $A$ only. An example is the classical cup
product $f,g \mapsto f \cup g$ of Hochschild cochains, resp.~the
induced graded commutative associative multiplication on the
Hochschild cohomology of associative algebras~\cite{gerstenhaber:AM63}.  Our
definition excludes some operations that are also ``natural'' in some
sense, such as the degree zero unary operation defined as the
projection $\pi_n : \HAP * \to \HAP n$, $n \geq 0$, because this
operation is not induced by any natural cochain map. A precise
definition of natural operations is given in
Section~\ref{natural}. Our aim is to {\em describe the homotopy type\/} of the
dg operad $\calB_\calP = \{\calB_\calP(n)\}_{n \geq 0}$ of all these
natural operations, see Problem~\ref{prob1} and its baby version
Problem~\ref{prob3}. The reward would be an ultimate understanding of
the structure of the cohomology of a given type of algebras.

Our original hope was that the homotopy type of $\calB_\calP$ would be
that of another Koszul quadratic operad $\calQ_\calP$ determined by
$\calP$ in an explicit and simple manner. Examples we had the in mind
were $\calP = \As$ for which probably $\calQ_\calP = \Gerst$, the
operad for Gerstenhaber algebras, and $\calP = \Lie$ for which
probably $\calQ_\calP = \Lie$, the operad for Lie algebras.
Calculations presented in this note however show that the homotopy
type of $\calB_\calP$ is in general more complicated, therefore the
property that makes the homotopy type of $\calB_\calP$ for $\calP =
\As$ or $\calP = \Lie$ so nice must be finer than just the Koszulity
of $\calP$. We have no idea what this property is.

We feel that our formulations are somehow unsatisfactory -- 
we would certainly prefer a concept that would not depend on a
``representation'' of the cohomology by a concrete cochain complex.
In an ideal world, we should be working with 
natural operations in an appropriate ``derived'' category in which 
the cohomology is the hom-functor. The
possibility of such a more conceptual approach for associative algebras and 
their Hochschild cohomology was indicated
by~\cite{kock-toen}, see also~\cite{hu,hu-kriz}. 

Another possibility could be to consider $H_\calP^*(A;A)$ as the
cohomology of the cotangent complex of a suitable suitably derived
stack of the variety of structure constants of $\calP$-algebras,
see~\cite{fontanine-kapranov:JAMS,fontanine-kapranov,markl:JPAA96},
and study automorphisms of the point of this stack representing the
algebra $A$. Our feeling is, however, that these fancier approaches,
despite their beauty and generality, are still not developed enough to
give {\em concrete\/} answers to {\em concrete\/} questions.

Let us explain the title of this note. In his famous
letter~\cite{deligne:letter}, P.~Deligne asked whether the
Gerstenhaber algebra structure on the Hochschild cohomology of an
associative algebra given by the cup product and the intrinsic bracket
is induced by a natural action of singular chains on the little discs
operad. There are several proofs of this so-called {\em Deligne
conjecture\/}
today~\cite{kaufmann,kontsevich-soibelman,mcclure-smith,mcclure-smith:JAMS03,%
tamarkin:deligne,hinich:FM03}. Assume one can prove that the
operad of all natural operations on the Hochschild complex (that is,
$\calB_\As$ in our notation) has the homotopy type of the operad for
Gerstenhaber algebras. The formality of the little discs
operad~\cite{tamarkin:formality} would then immediately imply the
Deligne conjecture by simple homological considerations.

In fact, most of the proofs of the Deligne conjecture we are aware
of~\cite{kaufmann,kontsevich-soibelman,mcclure-smith,mcclure-smith:JAMS03},
involve a conveniently chosen suboperad of $\calB_\Ass$ whose homotopy
type is detected by Fiedorowicz' recognition principle for
$E_2$-operads~\cite{fiedorowicz}. We will discuss these proofs in
Section~\ref{pisi_na_novem_notebooku}. Other proofs based on the
Etingof-Kazhdan (de)quantization were given
in~\cite{tamarkin:deligne,hinich:FM03}. Several attempts have
also been made to find a suitable filtration of the Fulton-MacPherson
compactification of the configuration space of points in the plane to
prove the conjecture~\cite{getzler-jones:preprint,voronov:99}. The
Deligne conjecture has surprising implications for the existence of
the deformation quantization of Poisson
manifolds~\cite{hinich:FM03,tamarkin:deligne}.

{\bf Acknowledgments.} I would like to express my thanks to
F.~Chapoton, E.~Getzler, V.~Hinich M.~Livernet,
P.~Somberg and A.A.~Voronov for many useful comments and remarks.  
My special thanks are due to D.~Tamarkin for inspiring discussions
during my stay at the Northwestern University in April~2004.

\section{Formulation of the problem}

In this section we state the problems sketched out in the
introduction more concretely and formulate also some conjectures.  
Let $\calB_\calP =
(\calB_\calP,\delta_\calP)$ be the dg-operad of all natural
multilinear operations on the cochain complex $C^*_\calP(A;A) =
(C^*_\calP(A;A),d_\calP)$.  The $n$-th component
$\calB_\calP(n)$ of $\calB_\calP$ is the space of all $n$-linear
natural operations $C^*_\calP(A;A)^{\ot n} \to \CAP*$ with the grading
induced by the grading of $C^*_\calP(A;A)$: $U \in \calB_\calP(n)$ has
degree $d$ if
\[
U(f_1,\ldots,f_n) \in \CPAA{m_1+\cdots +m_n + d},
\]
whenever $f_i \in \CPAA {m_i}$ for $1 \leq i \leq n$. In this case we
write $U \in \calB^d_\calP(n)$. Each $\calB_\calP(n)$ is equipped with
the degree $+1$ differential $\delta_\calP$ induced by the
differential $d_\calP$ of $C^*_\calP(A;A)$ in the usual way.

A precise definition of the operad $\calB_\calP$ is given in
Section~\ref{natural}.  Here we emphasize only that $\calB^d_\calP(n)
= 0$ for $d < 0$ and that $\calB_\calP(0) \not = 0$ for any nontrivial
$\calP$. The central problem of the paper reads:

\begin{problem}
\label{prob1}
Describe the homotopy type (in the non-abelian derived category) of
the dg operad $\calB_\calP$. In particular, calculate the cohomology
of $\calB_\calP$.
\end{problem}

A baby-version of this problem is Problem~\ref{prob3} of
Section~\ref{spaci_medvidek}. Closely related is:

\begin{problem}
\label{zase_prace_nad_hlavu}
Find a property characterizing operads $\calP$ for which $\calB_\calP$
is formal and has the homotopy type of some Koszul quadratic operad.
\end{problem}

We will see, in Example~\ref{ve_stredu_poletim_do_Minneapolis}, a
simple quadratic Koszul operad $\calD$ such that
\hbox{$H^*(\calB_\calD(0),\delta_\calD) \not= 0$}. This clearly means
that $\calB_\calD$ does not have the homotopy type of a quadratic Koszul
operad, therefore the property answering
Problem~\ref{zase_prace_nad_hlavu} must be stronger than Koszulness of
$\calP$.

Suppose that $\calP$ is the symmetrization of a non-$\Sigma$ operad
$\ucalP$~\cite[Remark~II.1.15]{markl-shnider-stasheff:book}. In this
case there exists a dg-suboperad $\calB_\ucalP$ of $\calB_\calP$
consisting of natural operations that preserve the
order of inputs of $\calP$-cochains.  For example, the classical cup product
$f\cup g \in C^1_\Ass(A;A) \cong \Lin(\otexp A2,A)$ of Hochschild
cochains $f,g \in C^0_\Ass(A;A) \cong \Lin(A,A)$ defined as 
\[
(f \cup g)(a \ot b) := f(a)\cdot g(b)\ 
\mbox { for } a \ot b \in A \ot A,
\]
with $\cdot$ denoting the associative multiplication of $A$, belongs to
$\calB_\uAss$, while the operation 
\[
U(f,g):= f(b)\cdot g(a)\ \mbox { for } a \ot b \in A \ot A,
\]
does not, see Definition~\ref{letel_jsem_do_Nymburka!!}
of Section~\ref{natural} for details.

Since $\calB_\ucalP(n)$ is a $\Sigma_n$-closed
subspace of $\calB_\calP(n)$, $n \geq 0$, $\calB_\ucalP$
is a usual, not only a non-$\Sigma$, operad.  We will see in
Example~\ref{zase_nebude_pocasi} that, surprisingly, the homotopy type of
$\calB_\ucalP$ in general {\em differs\/} from the homotopy type of
$\calB_\calP$. We therefore formulate:

\begin{problem}
\label{Pastviny}
Let $\ucalP$ a non-$\Sigma$ quadratic Koszul operad.  Describe the
homotopy type of the dg operad $\calB_\ucalP$. In particular,
calculate the cohomology of $\calB_\ucalP$.
\end{problem}

In Section~\ref{s6}(i) we give some indications that the operad
$\calB_\uAss$ has the homotopy type of the operad $\Gerst$ for
Gerstenhaber algebras, see~A.\ref{ma_prejit_fronta_v_sobotu} for a
definition of $\Gerst$.

One may consider also {\em strongly homotopy versions\/} of the above
problems.  Recall that a strongly homotopy $\calP$-algebra is,
by~\cite{markl:zebrulka}, an algebra over the minimal model $sh\calP$
of the operad $\calP$. Let us denote by $sh\calB_\calP =
\calB_{sh\calP}$ the dg-operad of natural operations on the cochain
complex $C^*_{sh\calP}(A,A)$ for the cohomology of a strongly homotopy
algebra $A$ with coefficients in itself. An example of this type of
operad is the ``minimal operad'' $M$ considered
in~\cite{kontsevich-soibelman}, which is a certain suboperad of
$\calB_{sh\Ass}$, see Section~\ref{pisi_na_novem_notebooku}(iii).

It is clear that there exists a canonical map $\calB_{\calP} \to
\calB_{sh\calP}$, but simple examples show that, again rather
surprisingly, this map is in general not a homotopy equivalence. Let
us formulate:

\begin{problem}
Describe the homotopy type of the dg-operad $\calB_{sh\calP}$ of
natural operations on the cohomology of strongly homotopy $\calP$-algebras.
\end{problem}

 Other problems formulated in
this paper are Problem~\ref{-} of Section~\ref{s2} and
Problems~\ref{prob3},\ref{+} of Section~\ref{s3}.

Let us finally formulate also some conjectures.
Although the operads $\calB_\calP$ and $\calB_{\calP^!}$ are
not isomorphic (see Section~\ref{natural}), computational evidences
together with an equivalence between the derived category of $\calP$
algebras and the derived category of $\calP^!$-algebras lead us to
believe in:

\begin{conjecture}
\label{zase}
The homotopy type of the operad $\calB_\calP$ is the same as the
homotopy type of $\calB_{\calP^!}$.
\end{conjecture}

The following two conjectures concern the homotopy type of
$\calB_\calP$ for $\calP = \Ass$ and $\calP = \Lie$.

\begin{conjecture}
The operad $\calB_\Ass$ has the homotopy type of the operad $\Gerst$ for
Gerstenhaber algebras.
\end{conjecture}

Some results which may be helpful in the proof of the above
conjecture are recalled in Section~\ref{pisi_na_novem_notebooku}.

\begin{conjecture}
\label{Kuratko}
The operad $\calB_\Lie$ has the homotopy type of the operad $\Lie$.
\end{conjecture}

According to a formality theorem~\cite[Proposition~3.4]{markl:zebrulka}, it
is enough to prove that 
\[
H^*(\calB_\Lie,\delta_\Lie) \cong \Lie.
\]
Since $H^0(\calB_\Lie,\delta_\Lie) \cong \Lie$ (see
Section~\ref{jsem_otrok}), Conjecture~\ref{Kuratko} is equivalent to
the acyclicity of $\calB_\Lie$ in positive degrees. Another conjecture,
Conjecture~\ref{ib}, is given in Section~\ref{s4}.

Let us finish this section with one exceptional example.
The trivial operad $\triv$ is a Koszul quadratic
self-dual operad. A $\triv$-algebra is a vector space $A$ with no
operations. Clearly $C^*_\triv(A;A)$ is just the space $\Lin(A,A)$ of
linear maps $f : A \to A$
concentrated in degree zero with trivial differential, thus
$H^*_\triv(A;A) = \Lin(A,A)$. It is also clear that
all natural operations on $\Lin(A,A)$ are the identity $\id_A \in
\Lin(A,A)$ considered as a degree zero constant, and 
iterated compositions 
\[
\Lin(A,A) \ni f_1,f_2,\ldots,f_n \mapsto f_1 \circ f_2 \circ
\cdots \circ f_n \in \Lin(A,A),\ n \geq 1.
\]
Therefore
\[
\calB_\triv \cong U\hskip -.2em  \As,
\] 
the operad for unital associative algebras. This example is
pathological in that the canonical element introduced in
Definition~\ref{canon} equals zero. Therefore, from now on all
quadratic Koszul operads in this note will be {\em nontrivial\/} in
the sense that $\calP \not= \triv$.

\section{The constants $\calB_\calP(0)$ -- soul without body}
\label{duneni_subwooferu}
\label{s2}

This section, as well as the rest of the paper, relies on terminology
and notation recalled in the Appendix.  The main result of this part
is Proposition~\ref{kozulka} which describes the dg-vector space
$\calB_\calP(0) = (\calB_\calP(0),\delta_\calP)$ of ``constants.''  It
is not difficult to see (compare also Example~\ref{Kocicka_zrcatko} of
Section~\ref{natural}) that
\[
\calB^{m-1}_\calP(0)  \cong \ss (\calP(m) \ot \calP^!(m))^{\Sigma_m}, \ m
\geq 1,
\]
with the action $\calB_\calP^{m-1}(0) \to \CAP{m-1}$ given as the composition
\begin{eqnarray}
\label{napsala_mi_Jituska!}
\lefteqn{ 
\ss(\calP(m) \ot \calP^!(m))^{\Sigma_m}
\stackrel{\cong}{\longrightarrow}
(\ss\calP(m) \ot \calP^!(m))^{\Sigma_m}
\stackrel{\ss \alpha \ot \id}{\longrightarrow}
(\ss\End_A(m) \ot \calP^!(m))^{\Sigma_m}} \hskip 1em
\\
\nonumber 
&&
\hskip 1em
\cong 
(\End_{\downarrow   A}(m) \ot \calP^!(m))^{\Sigma_m}
=
\left[
\Lin ((\downarrow \hskip -.2em A)^{\ot m},
\downarrow \hskip -.2em A) \ot \calP^!(m)\right]^{\Sigma_m}
=
\CAP{m-1}.
\end{eqnarray}

Since composition~(\ref{napsala_mi_Jituska!}) is monic for all
``generic'' $\calP$-algebras $A$, $(\calB^*_\calP(0),\delta_\calP)$ is
``morally'' the subcomplex of natural elements in $(\CAP *,d_\calP)$.
Before going further, we must recall the following general
construction.  Let $\calT$ be an operad. It is well-known that the
formula
\[
[f,g] := f \circ g - (-1)^{(m-1)(n-1)} g \circ f,
\]
where $f \circ g$ is, for $f \in \calT(m)$ and $g \in \calT(n)$,
defined by
\[
f \circ g := \sum_{1 \leq i \leq m} (-1)^{(n-1)(i-1)} f \circ_i g,
\]
makes the direct sum $\calT_* = \bigoplus_{m \geq
0}\calT_*$, with 
\[
\calT_{m-1} := \susp^{m-1} \calT(m) = \ss \calT(m),
\] 
a graded Lie algebra. Another standard fact is that each
element $\omega \in \calT_1 = \ss{}\calT(2)$ satisfying
$[\omega,\omega] = 0$ defines a degree $+1$ differential
$\delta_\omega : \calT_* \to \calT_{*+1}$ by
\[
\delta_\omega(t) := [t,\omega], \ \mbox { for } t \in \calT_*.
\]
It is helpful to observe that the condition $[\omega,\omega] =
0$ means the associativity:
\begin{equation}
\label{vcera_jsme_se_koupali_na_Hradistku}
\omega \circ_1 \omega = \omega \circ_2 \omega
\end{equation}
and that the differential $\delta_\omega$ in terms of
$\circ_i$-operations equals
\[
\delta_{\omega}(t) =  t \circ_1 \omega - t \circ_2 \omega
+ \cdots -(-1)^m  t \circ_m \omega + (-1)^{m}\omega \circ_1 t  - 
\omega \circ_2 t,
\ \mbox { for } t \in \calT(m).
\]

As proved in~\cite{laan-defo}, the graded Lie algebra structure
$(\calT_*,[-,-])$ descents to the space of coinvariants therefore it
induces, via the canonical isomorphism between invariants and
coinvariants, a Lie bracket, denoted again $[-,-]$, on the graded
vector space $\calT^\Sigma_* = \bigoplus_{m \geq 0}\calT^\Sigma_m$
with pieces
\[
\calT^\Sigma_{m-1} := \hskip .2em \uparrow^{m-1} 
(\calT(m) \ot \sgn_m)^{\Sigma_{m}} = \ss \calT(m)^{\Sigma_{m}}.
\]
As usual, an element $\phi \in \calT^\Sigma_1 
= \ss \calT(2)^{\Sigma_2}$ satisfying
$[\phi,\phi]=0$ induces a degree $+1$ differential $\delta^\Sigma_\phi:
\calT^\Sigma_* \to \calT^\Sigma_{*+1}$ by
\begin{equation}
\label{dnes_poslan_ib_do_MA}
\delta^\Sigma_\phi t := [\phi,t],\ \mbox { for } t \in \calT^\Sigma_*.
\end{equation}

In Proposition~\ref{kozulka} below we put $\calT := (\calP \ot
\calP^!)$ and define the differential~(\ref{dnes_poslan_ib_do_MA}) by
taking as $\phi$ the canonical element $\chi$ introduced in the
following definition in which $\#$ denotes the linear dual.

\begin{definition}
\label{canon}
Let $\calP$ be a quadratic Koszul operad. The {\em canonical
element\/} $\chi$ is the element of $\ss(\calP \ot
\calP^!)(2)^{\Sigma_2}$ corresponding, under the standard
identification 
\[
\ss(\calP \ot \calP^!)(2) 
\cong \hskip .2em \uparrow \hskip -.2em (\calP \ot \calP^\#)(2) 
\cong \hskip .2em \uparrow \hskip -.2em (\calP(2) \ot \calP(2)^\#) 
\cong \hskip .2em \uparrow \hskip -.2em\Lin(\calP(2), \calP(2)),
\]
to the suspension of the identity map $\uparrow \hskip -.2em
\id_{\calP(2)} \in \hskip .2em\uparrow \hskip -.2em\Lin(\calP(2),\calP(2))$.
\end{definition}
Observe that $\chi$ is symmetric,
\begin{equation}
\label{porad_to_knika}
\chi \tau = \chi\ \mbox { for
$\tau \in \Sigma_2$,}
\end{equation}
therefore indeed $\chi \in \ss(\calP \ot \calP^!)(2)^{\Sigma_2}$. The
condition $[\chi, \chi]=0$ is equivalent to the Jacobi
identity~(\ref{zase_jsem_se_opil}) for $\chi$ which follows
from~\cite[Corollary~2.2.9(b)]{ginzburg-kapranov:DMJ94}, see also the proof of
Proposition~\ref{lisasek}.

\begin{proposition}
\label{kozulka}
There is a natural isomorphism of cochain complexes
\[
(\calB^*_\calP(0),\delta_\calP) \cong ( (\opt)^\Sigma_*,\delta^\Sigma_\chi).
\]
\end{proposition}

If $\calP$ is the symmetrization of a \ns\ operad
$\ucalP$~\cite[Remark~II.1.15]{markl-shnider-stasheff:book}, then
there is a similar description of the chain complex
$(\calB^*_\ucalP(0),\delta_\ucalP)$ obtained as follows. 
The definition of the graded Lie algebra $(\calT_*, [-,-])$ given
above clearly makes sense also when $\calT$ is a non-$\Sigma$ operad. Observe
also that there exists the \ns\ quadratic dual $\ucalP^!$ of $\ucalP$
and that one may introduce the {\em \ns\ canonical element\/} $\uchi
\in \ss (\ucalP \ot \ucalP^!)(2)$ in exactly the same manner as its
symmetric version. The element $\uchi$ obviously satisfies the associativity
condition~(\ref{vcera_jsme_se_koupali_na_Hradistku}):
\[
\uchi \circ_1 \uchi = \uchi \circ_2 \uchi.
\]
Our non-$\Sigma$ version of Proposition~\ref{kozulka} reads:

\begin{proposition}
\label{Jituska}
Let $\calP$ be the symmetrization of a quadratic Koszul non-$\Sigma$
operad $\ucalP$.  Then
\[
(\calB^*_\ucalP(0),\delta_\ucalP) \cong ((\ucalP \ot \ucalP^!)_*,\delta_\uchi).
\]
\end{proposition}

Let us make a comment on the meaning of the cohomology
$H^*(\calB_\calP(0),\delta_\calP)$.  The natural morphism
\[
M: H^*(\calB_\calP(0),\delta_\calP) \to \HAP*
\]
induced by action~(\ref{napsala_mi_Jituska!}) is monic for any
``generic'' $\calP$-algebra $A$, therefore elements
$H^*(\calB_\calP(0),\delta_\calP)$ represent {\em natural generically
nontrivial homology classes\/} in the cohomology of
$\calP$-algebras. This leads one to believe that
$H^*(\calB_\calP(0),\delta_\calP)=0$ for all well-behaved operads,
since otherwise people would stumble upon nontrivial natural classes -- compare
the Casimir element in the cohomology of simple Lie
algebras. Example~\ref{ve_stredu_poletim_do_Minneapolis} however
contradict this reasonable assumption. We believe that
$H^*(\calB_\calP(0),\delta_\calP)$ is an important invariant of the
operad $\calP$ that deserves its own name:

\begin{definition}
We call the graded vector space $H^*(\calB_\calP(0),\delta_\calP)$ described
in Proposition~\ref{kozulka} the
{\em soul of the cohomology\/} of $\calP$-algebras.
\end{definition}

It is easy to prove that $H^0(\calB_\calP(0),\delta_\calP)$
is always trivial.

\begin{example}
\label{ctvrek_a_zase_je_lat}
{\rm\
Let us describe the complex calculating the soul
$H^*(\calB_\As(0),\delta_\As)$ of the Hochschild cohomology. Clearly
\[
(\calP \ot \calP^!)^\Sigma_{m-1} = (\As \ot \As)^\Sigma_{m-1} \cong
\ss(\As \ot \As)(m)^{\Sigma_m} \cong \ss \As(m),
\]
therefore the complex $((\As \ot \As)^\Sigma,\delta^\Sigma_\chi)$ has the
form
\begin{equation}
\label{to_jsem_zvedav_jestli_budu_moct_ten_notebook_zaplatit_z_grantu}
\bfk \stackrel{\delta^\Sigma_\chi}{\longrightarrow} \bfk[\Sigma_2] 
\stackrel{\delta^\Sigma_\chi}{\longrightarrow}  \bfk[\Sigma_3] 
\stackrel{\delta^\Sigma_\chi}{\longrightarrow} \bfk[\Sigma_4]
\stackrel{\delta^\Sigma_\chi}{\longrightarrow} \cdots
\end{equation}
Is also easy to describe the differential $\delta^\Sigma_\chi$; on a permutation
$\sigma \in \Sigma_m$ it acts as
\[
\delta^\Sigma_\chi(\sigma) := d_0(\sigma) -  d_1(\sigma) +  d_2(\sigma) -
\cdots + (-1)^{m+1}  d_{m+1}(\sigma) \in \bfk[\Sigma_{m+1}],   
\]
where $ d_0(\sigma) := \id \times \sigma$, $d_{m+1}(\sigma) := \sigma
\times \id$ and $d_i(\sigma) \in \Sigma_{m+1}$ is the permutation
obtained by doubling the $i$th input of $\sigma$. In
Theorem~\ref{krystal_mame} below we prove
that~(\ref{to_jsem_zvedav_jestli_budu_moct_ten_notebook_zaplatit_z_grantu})
is acyclic.

Since $\Ass$ is the symmetrization of the non-$\Sigma$ operad $\uAss$,
it makes sense to consider also the subcomplex
$(\calB^*_\uAss(0),\delta_\uAss)$ of $(\calB^*_\Ass(0),\delta_\Ass)$
described in Proposition~\ref{Jituska}.  This subcomplex is obviously
isomorphic to the acyclic complex
\begin{equation}
\label{Andulka}
\bfk \stackrel{d_0}{\longrightarrow} \bfk
\stackrel{d_1}{\longrightarrow}  \bfk
\stackrel{d_2}{\longrightarrow} \bfk
\stackrel{d_3}{\longrightarrow} \cdots
\end{equation}
in which $d_{2i} = \id_\bfk$ and $d_{2i+1} = 0$, $i \geq 0$. The
inclusion $(\calB^*_\uAss(0),\delta_\uAss) \hookrightarrow
(\calB^*_\Ass(0),\delta_\Ass)$ sends the generator $1 \in \bfk$ of the
$n$th piece of~(\ref{Andulka}) into the identity permutation
$\id_{n-1} \in \bfk[\Sigma_{n-1}]$ in
complex~(\ref{to_jsem_zvedav_jestli_budu_moct_ten_notebook_zaplatit_z_grantu}).
}
\end{example}

\begin{theorem}
\label{krystal_mame}
The soul $H^*(\calB_\As(0),\delta_\As)$ of the Hochschild cohomology
is trivial.
\end{theorem}

\noindent
{\bf Proof.}  We must prove
that~(\ref{to_jsem_zvedav_jestli_budu_moct_ten_notebook_zaplatit_z_grantu})
is an acyclic complex. The idea will be to show that it
decomposes into a direct sum of acyclic subcomplexes indexed by
primitive, in the sense introduced below, permutations.

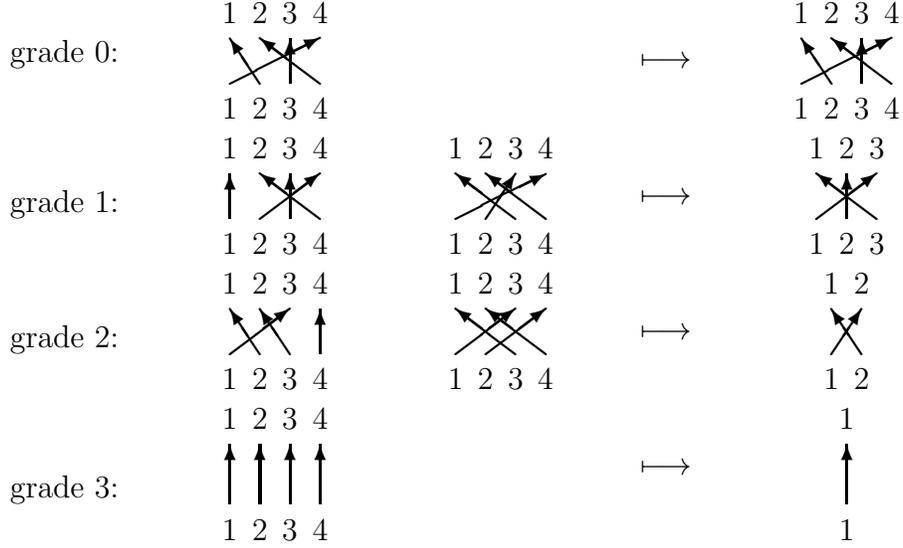
\begin{figure}
\begin{center}
\unitlength 2mm 
\thicklines
\begin{picture}(65,35)(4,5)
\put(21,35){\vector(2,1){6}}
\put(59,35){\vector(2,1){6}}
\put(23,35){\vector(-2,3){2}}
\put(61,35){\vector(-2,3){2}}
\put(27,35){\vector(-4,3){4}}
\put(65,35){\vector(-4,3){4}}
\put(25,35){\vector(0,1){3}}
\put(63,35){\vector(0,1){3}}
\put(21,26){\vector(0,1){3}}
\put(23,26){\vector(4,3){4}}
\put(27,26){\vector(-4,3){4}}
\put(25,26){\vector(0,1){3}}
\put(36,26){\vector(2,1){6}}
\put(40,26){\vector(-4,3){4}}
\put(42,26){\vector(-4,3){4}}
\put(38,26){\vector(2,3){2}}
\put(27,17){\vector(0,1){3}}
\put(36,17){\vector(4,3){4}}
\put(38,17){\vector(4,3){4}}
\put(42,17){\vector(-4,3){4}}
\put(40,17){\vector(-4,3){4}}
\put(21,17){\vector(4,3){4}}
\put(25,17){\vector(-2,3){2}}
\put(23,17){\vector(-2,3){2}}
\put(60,26){\vector(4,3){4}}
\put(64,26){\vector(-4,3){4}}
\put(62,26){\vector(0,1){3}}
\put(63,17){\vector(-2,3){2}}
\put(61,17){\vector(2,3){2}}
\put(21,7){\vector(0,1){4}}
\put(62,7){\vector(0,1){4}}
\put(23,7){\vector(0,1){4}}
\put(25,7){\vector(0,1){4}}
\put(27,7){\vector(0,1){4}}
\put(21,39){\makebox(0,0)[bc]{$1$}}
\put(59,39){\makebox(0,0)[bc]{$1$}}
\put(21,34){\makebox(0,0)[tc]{$1$}}
\put(59,34){\makebox(0,0)[tc]{$1$}}
\put(60,25){\makebox(0,0)[tc]{$1$}}
\put(60,30){\makebox(0,0)[bc]{$1$}}
\put(61,21){\makebox(0,0)[bc]{$1$}}
\put(61,16){\makebox(0,0)[tc]{$1$}}
\put(63,16){\makebox(0,0)[tc]{$2$}}
\put(21,30){\makebox(0,0)[bc]{$1$}}
\put(21,25){\makebox(0,0)[tc]{$1$}}
\put(21,21){\makebox(0,0)[bc]{$1$}}
\put(21,16){\makebox(0,0)[tc]{$1$}}
\put(21,12){\makebox(0,0)[bc]{$1$}}
\put(62,12){\makebox(0,0)[bc]{$1$}}
\put(21,6){\makebox(0,0)[tc]{$1$}}
\put(62,6){\makebox(0,0)[tc]{$1$}}
\put(23,39){\makebox(0,0)[bc]{$2$}}
\put(61,39){\makebox(0,0)[bc]{$2$}}
\put(23,30){\makebox(0,0)[bc]{$2$}}
\put(23,21){\makebox(0,0)[bc]{$2$}}
\put(23,34){\makebox(0,0)[tc]{$2$}}
\put(61,34){\makebox(0,0)[tc]{$2$}}
\put(62,25){\makebox(0,0)[tc]{$2$}}
\put(62,30){\makebox(0,0)[bc]{$2$}}
\put(63,21){\makebox(0,0)[bc]{$2$}}
\put(23,25){\makebox(0,0)[tc]{$2$}}
\put(23,16){\makebox(0,0)[tc]{$2$}}
\put(23,12){\makebox(0,0)[bc]{$2$}}
\put(23,6){\makebox(0,0)[ct]{$2$}}
\put(38,21){\makebox(0,0)[bc]{$2$}}
\put(38,25){\makebox(0,0)[tc]{$2$}}
\put(38,30){\makebox(0,0)[bc]{$2$}}
\put(38,16){\makebox(0,0)[tc]{$2$}}
\put(25,39){\makebox(0,0)[bc]{$3$}}
\put(63,39){\makebox(0,0)[bc]{$3$}}
\put(25,34){\makebox(0,0)[tc]{$3$}}
\put(63,34){\makebox(0,0)[tc]{$3$}}
\put(64,25){\makebox(0,0)[tc]{$3$}}
\put(64,30){\makebox(0,0)[bc]{$3$}}
\put(25,30){\makebox(0,0)[bc]{$3$}}
\put(25,25){\makebox(0,0)[tc]{$3$}}
\put(40,30){\makebox(0,0)[bc]{$3$}}
\put(40,25){\makebox(0,0)[tc]{$3$}}
\put(40,21){\makebox(0,0)[bc]{$3$}}
\put(40,16){\makebox(0,0)[tc]{$3$}}
\put(25,12){\makebox(0,0)[bc]{$3$}}
\put(25,6){\makebox(0,0)[ct]{$3$}}
\put(25,21){\makebox(0,0)[bc]{$3$}}
\put(25,16){\makebox(0,0)[tc]{$3$}}
\put(27,39){\makebox(0,0)[bc]{$4$}}
\put(65,39){\makebox(0,0)[bc]{$4$}}
\put(27,34){\makebox(0,0)[tc]{$4$}}
\put(65,34){\makebox(0,0)[tc]{$4$}}
\put(27,30){\makebox(0,0)[bc]{$4$}}
\put(27,25){\makebox(0,0)[tc]{$4$}}
\put(42,30){\makebox(0,0)[bc]{$4$}}
\put(42,25){\makebox(0,0)[tc]{$4$}}
\put(42,21){\makebox(0,0)[bc]{$4$}}
\put(42,16){\makebox(0,0)[tc]{$4$}}
\put(27,12){\makebox(0,0)[bc]{$4$}}
\put(27,6){\makebox(0,0)[ct]{$4$}}
\put(27,16){\makebox(0,0)[tc]{$4$}}
\put(27,21){\makebox(0,0)[bc]{$4$}}
\put(36,30){\makebox(0,0)[bc]{$1$}}
\put(36,25){\makebox(0,0)[tc]{$1$}}
\put(36,21){\makebox(0,0)[bc]{$1$}}
\put(36,16){\makebox(0,0)[tc]{$1$}}
\put(50,36){\makebox(0,0)[cb]{$\longmapsto$}}
\put(50,27){\makebox(0,0)[cb]{$\longmapsto$}}
\put(50,18){\makebox(0,0)[cb]{$\longmapsto$}}
\put(50,9){\makebox(0,0)[cb]{$\longmapsto$}}
\put(10,36){\makebox(0,0)[cb]{grade $0$:}}
\put(10,27){\makebox(0,0)[cc]{grade $1$:}}
\put(10,18){\makebox(0,0)[cc]{grade $2$:}}
\put(10,8){\makebox(0,0)[cc]{grade $3$:}}
\end{picture}
\end{center}
\caption{\label{Jitka_prisla}%
Left: Examples of elements of $\Sigma_4$ of grade $0$ (first line),
grade $1$ (second line), grade $2$ (third line) and grade $3$ (fourth
line). Right: the corresponding primitive elements.}
\end{figure}

We define first, for each $\sigma \in \Sigma_n$, a natural number
$g(\sigma)$, $0 \leq g(\sigma) \leq n$, which we call the {\em
grade\/} of $\sigma$, as follows.  The grade of the unit $\id_n \in
\Sigma_n$ is $n-1$, $g(\id_n) := n-1$.  For $\sigma \not= \id_n$, let
\begin{eqnarray*}
a(\sigma) &:=& \max\{i;\ \sigma = \id_i \times \tau \mbox { for some }
\tau \in \Sigma_{n-i}\}, \mbox { and }
\\
c(\sigma) &:=& \max\{j;\ \sigma = \lambda \times \id_j \mbox { for some }
\lambda \in \Sigma_{n-i}\}.
\end{eqnarray*}
There clearly exists a unique $\omega = \omega(\sigma) \in \Sigma_{n
  - a(\sigma) - c(\sigma)}$ such that
$\sigma = \id_{a(\sigma)} \times \omega(\sigma) \times
\id_{c(\sigma)}$.
Let, finally, $b(\sigma)$ be the number of ``doubled strings'' in
$\omega(\sigma)$, 
\[
b(\sigma) := \{1 \leq s < k;\ \omega(s+1) = \omega(s) +1\}.
\]
The grade of $\sigma$ is then defined by 
\[
g(\sigma) := a(\sigma) + b(\sigma) + c(\sigma),
\]
see Figure~\ref{Jitka_prisla} for examples. Observe that the
differential $\delta^\Sigma_\chi$
of~(\ref{to_jsem_zvedav_jestli_budu_moct_ten_notebook_zaplatit_z_grantu})
raises the grade by $+1$.

Let us call $\chi \in \Sigma_k$, $k \geq 1$, {\em primitive\/} if
$g(\sigma) = 0$. Observe that, according to our definitions, $\id_n
\in \Sigma_n$ is primitive if and only if $n=1$.  For each $\sigma \in
\Sigma_n$, $\sigma \not= \id_n$, we define a unique primitive $\kappa
= \kappa(\sigma) \in \Sigma_k$, $k = n - g(\sigma)$, by contracting all
``multiple strings'' of $\omega(\sigma)$ into ``simple'' ones, see
Figure~\ref{Jitka_prisla}. We put $\chi(\id_n):= \id_1$.

For a primitive $\chi$, let $P^*(\chi)$ be the graded subspace
of~(\ref{to_jsem_zvedav_jestli_budu_moct_ten_notebook_zaplatit_z_grantu})
spanned by all permutations $\sigma$ with $\chi = \kappa(\sigma)$. The
following statements can be easily verified:

(i) Each $P^*(\chi)$ is a subcomplex
of~(\ref{to_jsem_zvedav_jestli_budu_moct_ten_notebook_zaplatit_z_grantu}).

(ii)
Complex~(\ref{to_jsem_zvedav_jestli_budu_moct_ten_notebook_zaplatit_z_grantu})
decomposes as the summation $\bigoplus_\chi P^*(\chi)$
over all primitive permutations~$\chi$.

(iii) For each primitive $\chi \in \Sigma_n$,
\[
P^*(\chi) \cong P^*(\id_1) \ot \cdots \ot P^*(\id_1)\ \mbox { ($n+2$ times).}
\]
The proof is finished by observing that $P^*(\id_1)$ is isomorphic
to the acyclic complex~(\ref{Andulka}) and applying the K\"unneth formula.
\qed

\begin{example}
{\rm\ 
In this example we describe the soul of the Chevalley-Eilenberg
cohomology of Lie algebras which is, due to the obvious self-duality
of Proposition~\ref{kozulka}, the same as the soul of the Harrison
cohomology of commutative associative algebras. In both cases
\[
(\calP \ot \calP^!)^\Sigma_{m-1} = (\Com \ot \Lie)^\Sigma_{m-1} \cong
\ss \Lie(m)^{\Sigma_m} = \hskip .4em \uparrow^{m-1} \hskip -.2em
({\Lie(m) \otimes \sgn_m})^{\Sigma_m} \cong \bfk
\] 
(see~\cite{klya}) and one may identify $((\Com \ot
\Lie)^\Sigma,\delta^\Sigma_\chi)$ with
the acyclic complex~(\ref{Andulka}). Therefore the
souls of both the Chevalley-Eilenberg cohomology and the
Harrison cohomology are trivial.
}
\end{example}

\begin{example}
\label{ve_stredu_poletim_do_Minneapolis}
{\rm\ 
This example presents a Koszul quadratic operad with a nontrivial soul.
Let $\calD := \Ass * \Ass$ be the free product of two copies of
the associative operad. Operad $\calD$ is a
Koszul quadratic operad, 
whose quadratic dual ${\calD^!}$ equals the coproduct
$\Ass\vee\Ass$ defined by
\[
(\Ass\vee\Ass{})(m) := \cases{\bfk}{if $m=1$
                                  and}{\Ass(m)\oplus\Ass(m)}{if $m \geq 2$.}
\]
\begin{figure}[t]
\begin{center}
{
\unitlength=1.000000pt
\begin{picture}(230.00,175.00)(73.00,-10.00)
\put(0.00,120.00){\makebox(0.00,0.00){$d_2$}}
\put(0.00,70.00){\makebox(0.00,0.00){$d_2$}}
\put(0.00,20.00){\makebox(0.00,0.00){$d_2$}}
\put(270.00,10.00){\makebox(0.00,0.00)[l]{$d_1$}}
\put(160.00,10.00){\makebox(0.00,0.00)[l]{$d_1$}}
\put(50.00,10.00){\makebox(0.00,0.00)[l]{$d_1$}}
\put(380.00,0.00){\makebox(0.00,0.00){$\cdots$}}
\put(10.00,170.00){\makebox(0.00,0.00){$\vdots$}}
\put(330.00,0.00){\makebox(0.00,0.00){$\Lin(A^{\otimes 4},A)$}}
\put(10.00,150.00){\makebox(0.00,0.00){$\Lin(A^{\otimes 4},A)$}}
\put(220.00,0.00){\makebox(0.00,0.00){$\Lin(A^{\otimes 3},A)$}}
\put(10.00,100.00){\makebox(0.00,0.00){$\Lin(A^{\otimes 2},A)$}}
\put(110.00,0.00){\makebox(0.00,0.00){$\Lin(A^{\otimes 2},A)$}}
\put(10.00,50.00){\makebox(0.00,0.00){$\Lin(A^{\otimes 2},A)$}}
\put(10.00,0.00){\makebox(0.00,0.00){$\Lin(A,A)$}}
\put(10.00,110.00){\vector(0,1){30.00}}
\put(10.00,60.00){\vector(0,1){30.00}}
\put(10.00,10.00){\vector(0,1){30.00}}
\put(260.00,0.00){\vector(1,0){30.00}}
\put(150.00,0.00){\vector(1,0){30.00}}
\put(40.00,0.00){\vector(1,0){30.00}}
\end{picture}}
\end{center}
\caption{\label{jsem_nachlazeny}%
The ``meager'' bicomplex 
describing the cohomology of $\calD$-algebras.}
\end{figure} 

\noindent 
Obviously,
$\calD$-algebras are triples $A =
(A,\mu_1,\mu_2)$ consisting of two independent associative
multiplications 
$\mu_1,\mu_2 : A \ot A \to A$.
The cohomology of these algebras is the cohomology of the total
complex $(\CAD *,d_\calD)$ of the ``meager'' bicomplex
in Figure~\ref{jsem_nachlazeny}. The
horizontal line is the Hochschild cochain complex of the associative algebra
$A_1 := (A,\mu_1)$, the
vertical line the Hochschild complex of $A_2 := (A,\mu_2)$.

Let $e$ denote the identity $\id \in \Lin(A,A)$ considered as a natural
element of $\CAD 0$. Clearly
\[
d_\calD(d_1 e) = d_\calD(d_2 e) = 0
\]
therefore both $d_1 e$ and $d_2 e$ are natural cochains in $\CAD 1$
thus representing $\delta_\calD$-cochains in
$\calB_\calD^1(0)$. The equality
\[
d_\calD e = d_1 e + d_2 e
\]
implies that $d_1 e + d_2 e$ is $\delta_\calD$-homologous to zero in
$\calB^1_\calD(0)$. We conclude that
\[
H^1(\calB_\calD(0),\delta_\calD) \cong \Span([d_1 e]) \cong \Span([d_2 e]) \cong \bfk.
\] 
}
\end{example}

We saw that the souls of the Hochschild ($\calP = \As)$,
Chevalley-Eilenberg ($\calP = \Lie$) and Harrison ($\calP = \As$)
cohomologies were trivial, while the soul of the cohomology for
$\calD$-algebras analyzed in
Example~\ref{ve_stredu_poletim_do_Minneapolis} was not.  This leads us
to formulate:

\begin{problem}
\label{-}
Which property of a quadratic Koszul operad $\calP$ implies the
triviality of the soul $H^*((\calP \ot \calP^!)^\Sigma,\delta^\Sigma_\chi)$
of the $\calP$-cohomology?
\end{problem}

\begin{example}
\label{zase_nebude_pocasi}
{\rm\ In this example we describe a non-$\Sigma$ quadratic Koszul
operad with the property that $(\calB^*_\ucalP(0),\delta_\ucalP)$ is
acyclic but the soul $(\calB^*_\calP(0),\delta_\calP)$ is not. This
shows that the homotopy type of $\calB_\calP$ is in general different
from the homotopy type of $\calB_\ucalP$.
 
Let $\uMag$ be the free non-$\Sigma$ operad generated by one bilinear
operation, $\uMag : = \underline{\Gamma}( \underline{\mu})$, and $\Mag$ its
symmetrization. The corresponding cochain complex $C^*_\Mag (A;A)$ is
the truncation
\[
\Lin(A,A) \stackrel{d}{\longrightarrow} \Lin (\otexp A2,A)
\]  
of the Hochschild complex. The complex $(\calB^*_\Mag(0),\delta_\Mag)$
defining the soul of $\Mag$ is the truncation
\[
\bfk \stackrel{\delta^\Sigma_\chi}{\longrightarrow} \bfk[\Sigma_2]
\]
of~(\ref{to_jsem_zvedav_jestli_budu_moct_ten_notebook_zaplatit_z_grantu}),
and is manifestly not acyclic. On the other hand,
$(\calB^*_\uMag(0),\delta_\uMag)$ is acyclic, isomorphic to the truncation
$
\bfk \stackrel{d_0}{\longrightarrow} \bfk
$
of~(\ref{Andulka}). We conclude that
$H^*(\calB^*_\uMag(0),\delta_\uMag) = 0$ while
\[
H^*(\calB^*_\Mag(0),\delta_\Mag) = H^1(\calB^*_\Mag(0),\delta_\Mag)
\cong \bfk.
\]
}
\end{example}

\section{Homotopy type  of $\calB(1)$ -- surprises continue}
\label{spaci_medvidek}
\label{s3}

In this section we study, as a next step toward the understanding of
$\calB_\calP$, the homotopy type of the associative dg-algebra
$\calB_\calP(1) = (\calB^*_\calP(1),\delta_\calP)$.  Since the operad
$\calP^!$ is a module, in the sense of~\cite{markl:zebrulka}, over itself,
it makes sense to consider the space $\Endom_{\calP^!}({\calP^!})$ of
all ${\calP^!}$-module endomorphisms $\alpha : {\calP^!} \to
{\calP^!}$. Very crucially,
\begin{equation}
\label{volal_jsem_Martinovi}
\Endom_{\calP^!}({\calP^!}) \cong \bfk,
\end{equation}
because each $\alpha \in \Endom_{\calP^!}({\calP^!})$ is uniquely
determined by the value $\alpha_1(1) \in {\calP^!}(1)\cong \bfk$ and,
conversely, for each $\varphi \in \bfk$ the homomorphism $\alpha :=
\varphi \cdot \id_{\calP^!}$ is such that $\alpha_1(1) = \varphi$.

\begin{proposition}
\label{modules}
There is a {\em canonical\/} identification of associative unital algebras
\begin{equation}
\label{8}
H^0(\calB^*_\calP(1),\delta_\calP) 
\cong  \Endom_{\calP^!}({\calP^!}) \cong \bfk.
\end{equation}
\end{proposition}

\noindent 
{\bf Proof.}
Since there are no elements in negative degrees, 
\[
H^0(\calB_\calP(1)) = \Ker\left( \delta : \calB^0_\calP(1) \to
\calB^1_\calP(1)\right).
\]
By definition, elements of the kernel  $\Ker(\delta)$ are natural chain
maps
\[
\{\varphi_m : C^m_\calP(A;A) \to  C^m_\calP(A;A)\}_{m \geq 0}.
\]
As explained in Example~\ref{Deda_Mraz}, the naturality of $\varphi_m$
means that it is induced by a $\Sigma_{m+1}$-equivariant map
$\alpha_{m+1} : {\calP^!}(m+1) \to {\calP^!}(m+1)$ .  It is easy to
verify that the collection $\{\alpha_m\}$ determines a chain map if
and only if it assembles into a ${\calP^!}$-module endomorphism
$\alpha : {\calP^!} \to {\calP^!}$.
\qed

The following example shows that the dg-associative algebra
$\calB_\calP(1)$ might in general have nontrivial cohomology in
positive degrees.

\begin{example}
\label{sym}
{\rm Let $\Sym$ be the operad describing algebras with one commutative
bilinear multiplication and no axioms. Explicitly, $\Sym$ is the free
operad generated by the trivial representation of $\Sigma_2$ placed in
arity two.  It is a Koszul quadratic operad whose quadratic dual
$\Sym^!$ is given by $\Sym^!(1) = \bfk$, $\Sym^!(2) = \sgn_2$ (the
signum representation of $\Sigma_2$), and $\Sym^!(m) = 0$ for $m \geq
3$.

The cohomology of a $\Sym$-algebra $A = (A,\ \cdot\ )$ is the
cohomology of the two-term complex (which should be interpreted as a
truncation of the Harrison complex)
\[
\Lin(A,A) \stackrel{d}{\longrightarrow} \Lin(S^2 A,A),
\]
where $S^2 A$ is the second symmetric power of $A$.
The differential $d$ is given by the usual formula
\[
(d \phi)(a,b) := a \cdot \phi(b) - \phi(a\cdot b) + \phi(b) \cdot a,
\]
for $\phi \in \Lin(A,A)$ and $a,b \in A$.

We are going to describe the dg-algebra $\calB^*_\Sym(1)$.
Let $\alpha$ be the projection of $\Lin(A,A) \oplus \Lin(S^2 A,A)$
onto $\Lin(A,A)$ and $\beta$ the
projection onto $\Lin(S^2 A,A)$. Let $u$ and $v$ be degree $+1$
operations given by 
\[
u(\phi)(a,b) := a \cdot \phi(b) + \phi(a) \cdot b
\mbox { and }
v(\phi)(a,b) := \phi(a\cdot b),
\] 
for $\phi \in \Lin(A,A)$ and $a,b \in A$. Then clearly 
$\calB^0_\Sym(1)$ is the semisimple algebra $\bfk \oplus \bfk$ spanned by
$\alpha$ and $\beta$, 
and the space $\calB^1_\Sym(1)$ is two-dimensional, spanned by
$u$ and $v$. The higher $\calB^i_\Sym(1)$ are, for $i \geq 2$, trivial. To
describe the multiplication in $\calB^*_\Sym(1)$, it is enough to specify
how $\calB^0_\Sym(1)$ acts on $\calB^1_\Sym(1)$. This action is given by
\[
\alpha b = 0 = b\beta\ \mbox { and } b \alpha = b = \beta b, \ \mbox {
for $b \in \calB^1_\Sym(1)$.} 
\]
The differential $\delta_\Sym$ of $\calB^*_\Sym(1)$ acts by
\[
\delta \alpha = -\delta \beta = u - v,\
\delta u = \delta v = 0.
\]
The cohomology of $(\calB^*_\Sym(1),\delta_\Sym)$ can be easily calculated,
\[
H^*(\calB^*_\Sym(1)) \cong \bfk \oplus W,
\]
where $W$ is the vector space spanned by the class
$[u]$. We leave as a simple exercise to prove that there
exist a quasi-isomorphism $H^*(\calB^*_\Sym(1)) \to \calB^*_\Sym(1)$. 
The dg-associative algebra $\calB_\Sym(1)$ is therefore {\em formal\/}.
}
\end{example}

Here is a baby version of Problem~\ref{prob1}:

\begin{problem}
\label{prob3}
Describe the homotopy type of the unital differential graded
associative algebra $\calB_\calP(1) =
(\calB^*_\calP(1),\delta_\calP)$.  In particular, calculate the
cohomology of $\calB_\calP(1)$.
\end{problem}

We expect that the homotopy type of $\calB_\calP(1)$ is that of $\bfk$
for all ``reasonable'' operads, though we do not know what
``reasonable'' means -- the operad $\Sym$ of Example~\ref{sym} seems
reasonable enough, yet the homotopy type of $\calB_\Sym(1)$ is
nontrivial.
Let us close this section by formulating:

\begin{problem}
\label{+}
Which property of the operad $\calP$ implies that the
dg associative algebra $(\calB^*_\calP(1),\delta_\calP)$ has the
homotopy type of $\bfk$?
\end{problem}

\section{The operad $H^0(\calB_\calP)$ and the intrinsic bracket}
\label{jsem_otrok}
\label{s4}

It is well-known~\cite{markl:ib,schlessinger-stasheff:JPAA85} that the
chain complex $C^*_\calP(A;A)$ always carries a natural dg Lie algebra
structure given by the {\em intrinsic bracket\/}.  The easiest way to
construct such a bracket is to identify $C^*_\calP(A;A)$ with the dg
Lie algebra ${\it Coder\/}^*({\mathcal F}^c_{\calP^!}(\downarrow \hskip
-.2em A))$ of coderivations of the cofree nilpotent
$\calP^!$-coalgebra cogenerated by the desuspension $\downarrow \hskip
-.2em A$ as it was done
in~\cite[Definition~II.3.99]{markl-shnider-stasheff:book}.  In this
way we obtain a natural homomorphism
\begin{equation}
\label{lepidlo}
I: (\Lie,0) \to (\calB_\calP,\delta_\calP)
\end{equation}
of dg operads. If $\calP$ is the symmetrization of a non-$\Sigma$ operad
$\ucalP$, then $\Im(I) \subset \calB_\ucalP$, therefore the map $I$
of~(\ref{lepidlo}) factorizes through the natural inclusion
$(\calB_\ucalP,\delta_\ucalP) \hookrightarrow
(\calB_\calP,\delta_\calP)$. Computational evidences lead us to:

\begin{conjecture}
\label{ib}
The natural homomorphism $I: (\Lie,0) \to (\calB_\calP,\delta_\calP)$ 
induces an isomorphism of operads
\[
 H^0(\calB_\calP) \cong \Lie, 
\]
for an arbitrary nontrivial quadratic Koszul $\calP$.
\end{conjecture}

We were able to verify Conjecture~\ref{ib} for $\calP = \Lie$, that is,
to prove
\begin{equation}
\label{naspet_od_MF}
H^0(\calB_\Lie) \cong \Lie.
\end{equation}
This isomorphism turned out to be related to a certain
characterization of free Lie algebras inside free pre-Lie
algebras. More precisely, let $\freepreLie(X)$ denote the free pre-Lie
algebra generated by a set $X$. The commutator of the pre-Lie product
makes $\freepreLie(X)$ a Lie algebra. Let $\freeLie(X) \subset
\freepreLie(X)$ be the Lie algebra generated by $X$ in
$\freepreLie(X)$. It is not hard to see that $\freeLie(X)$ is in fact
isomorphic to the free Lie algebra on $X$, see also~\cite{dzhu-lof:HHA02}.
Then~(\ref{naspet_od_MF}) is implied by a very explicit
characterization of the subspace $\freeLie(X)$ of $\freepreLie(X)$.

Similarly, the conjectured isomorphism 
$H^0(\calB_\uAss) \cong \Lie$ can be
translated into a certain characterization of free Lie algebras inside free
brace algebras.
We were also able to prove that, for an arbitrary quadratic Koszul operad
\begin{equation}
\label{co_to_ten_Jim_ji?}
H^0(\calB_\calP(2)) \cong \sgn_2,
\end{equation}
the signum representation of $\Sigma_2$, by describing
$H^0(\calB_\calP(2))$ in terms of suitably defined {\em pairings\/}
$\calP^! \ot \calP^! \to \calP^!$.

Let us close this section by a couple of remarks which will be useful
in the proof of Proposition~\ref{lisasek}. As we recalled at the
beginning of this section, there is a canonical isomorphism
$C^*_\calP(A;A) \cong {\it Coder\/}^*({\mathcal F}^c_{\calP^!}(\downarrow 
\hskip -.2em A))$. It is well-known that coderivations of a cofree
nilpotent algebra form a natural pre-Lie
algebra~\cite[II.3.9]{markl-shnider-stasheff:book}, therefore one has
a natural homomorphism of non-dg operads
\begin{equation}
\label{nasadil_jsem_retez}
{\it pre}I : \preLie \to \calB_\calP.
\end{equation}
The map~(\ref{lepidlo}) is then the composition
\[
\Lie \longrightarrow \preLie \stackrel{{\it pre}I}{\longrightarrow}
\calB_\calP 
\]
of ${\it pre}I$ with the anti-symmetrization map $\Lie \to \preLie$. 

\section{The cup products}
\label{cup}
\label{s5}

The central statement of this section is Theorem~\ref{zitra_prednaska}
that claims that the suspension $\ss(\calP \ot \calP^!)$
(see~A.\ref{odpalil_jsem_notebook}) of the operad $(\calP \ot
\calP^!)$ acts on $\CAP*$, and Theorem~\ref{poslouchal_jsem_Kryla}
that characterizes which elements of $\ss(\calP \ot \calP^!)$ act via
chain maps.  Observe that the operad $\ss(\calP \otimes \calP^!)$ need
not be quadratic even when $\calP$ is.

\begin{theorem}
\label{zitra_prednaska}
There is a canonical action of the operad $\ss(\opt)$
on  the graded vector space $\CAP*$, via natural operations.
This action can be interpreted as an inclusion of non-differential
graded operads
\begin{equation}
\label{dnes_jsem_mel_prednasku_ve_Filadelfii}
\ccc: \ss(\calP \otimes \calP^!) \hookrightarrow \calB_\calP.
\end{equation}
\end{theorem}

\noindent 
{\bf Proof.}  The proof relies on the notation introduced/recalled
in~A.\ref{Ferda} and~A.\ref{zase_do_Prahy}.  
The ``tautological'' action of the endomorphism
operad $\End_A$ on $A$ tensored with the action of $\calP^!$ on itself
makes the graded vector space $\tildeCAP* = \bigoplus_{m \geq
0}\tildeCAP m$ an $\ss(\End_A \ot \calP^!)$-algebra. It is
straightforward to prove that this action induces, via
\begin{equation}
\label{zitra_letim_do_Raleigh}
t(f_1,\ldots,f_n) := {\it Aver\/}\left(\rule{0em}{1em}
t(\iota(f_1),\ldots,\iota(f_n))\right),
\end{equation}
for $t \in \ss(\End_A \ot\calP^!)(n)$ and $f_1,\ldots,f_n \in \CAP*$,
an action of $\ss(\End_A \ot \calP^!)$ on the graded vector space
$\CAP*$.  Suppose that $A$ is a $\calP$-algebra, with the structure
given by $\alpha : \calP \to \End_A$.
Action~(\ref{dnes_jsem_mel_prednasku_ve_Filadelfii}) is obtained by
composing action~(\ref{zitra_letim_do_Raleigh}) with the homomorphism
$\ss(\alpha \ot \id) : \ss(\calP \ot \calP^!) \to \ss(\End_A \ot
\calP^!)$.  An alternative description
of~(\ref{dnes_jsem_mel_prednasku_ve_Filadelfii}) is given in
Example~\ref{fit} of Section~\ref{natural}.%
\qed

We use inclusion~(\ref{dnes_jsem_mel_prednasku_ve_Filadelfii})
to view $\ss(\calP \ot \calP^!)$ 
as a suboperad of $\calB_\calP$.
Elements of $\ss(\calP \ot \calP^!)$ need not be $\delta_\calP$-closed 
in $\calB_\calP$;
let $\calZ_\calP \subset \ss(\calP \otimes \calP^!)$ 
denote the suboperad of $\delta_\calP$-cocycles. 
In Proposition~\ref{poslouchal_jsem_Kryla}, which describes
$\calZ_\calP$ explicitly, we use
the canonical element $\chi$ introduced in Definition~\ref{canon}.

\begin{theorem}
\label{poslouchal_jsem_Kryla}
The suboperad $\calZ_\calP$ of $\delta_\calP$-closed elements in
$\ss(\calP \ot \calP^!)$ is characterized as follows: $t \in \ss(\calP
\otimes \calP^!)(n)$ belongs to $\calZ_\calP(n)$ if and only if
\begin{equation}
\label{s_ni_jsem_se_natrapil}
\chi \circ_2 t + t \circ_1 \chi +   (t \circ_2 \chi) (12) +
 (t \circ_3 \chi) (123)  + \cdots +  (t\circ_n \chi) (123 \cdots n) = 0,
\end{equation}
where $(123\cdots k) \in \Sigma_{n+1}$ is the cycle
\[
\left(
\begin{array}{cccccccc}
1&2&3&\cdots&k&k+1&\cdots&n+1
\\
2&3&4&\cdots&1&k+1&\cdots&n+1
\end{array}
\right).
\]
\end{theorem}

The {\bf proof} is a completely straightforward calculation. We
recommend as an exercise to verify that solutions
of~(\ref{s_ni_jsem_se_natrapil}) are indeed closed under operadic
composition.  The meaning of equation~(\ref{s_ni_jsem_se_natrapil})
should be clear from Figure~\ref{jj}. The importance of the operad
$\calZ_\calP$ is explained by

\begin{figure}[t]
\begin{center}
{
\unitlength=1.000000pt
\begin{picture}(330.00,55.00)(0.00,10.00)
\thicklines
\put(330.00,40.00){\makebox(0.00,0.00){$=\ 0$}}
\put(240.00,40.00){\makebox(0.00,0.00){$+$}}
\put(160.00,40.00){\makebox(0.00,0.00){$+$}}
\put(70.00,40.00){\makebox(0.00,0.00){$+$}}
\put(280.00,30.00){\line(0,-1){10.00}}
\put(270.00,30.00){\line(0,-1){10.00}}
\put(280.00,20.00){\line(1,-2){10.00}}
\put(270.00,20.00){\line(1,-2){10.00}}
\put(-5,0){
\put(310.00,20.00){\makebox(0.00,0.00){$\chi$}}
\put(300.00,20.00){\makebox(0.00,0.00){$\bullet$}}
\put(300.00,20.00){\line(-3,-1){29.5}}
\put(310.00,10.00){\line(0,-1){10.00}}
\put(270.00,10.00){\line(0,-1){10.00}}
\put(300.00,20.00){\line(1,-1){10.00}}
}

\put(280.00,40.00){\makebox(0.00,0.00){$t$}}
\put(290.00,30.00){\line(1,-2){5.00}}
\put(280.00,50.00){\line(0,1){20.00}}
\put(280.00,50.00){\line(-1,-1){20.00}}
\put(300.00,30.00){\line(-1,1){20.00}}
\put(260.00,30.00){\line(1,0){40.00}}
%
\put(203.00,20.00){\makebox(0.00,0.00)[lb]{$\chi$}}
\put(200.00,40.00){\makebox(0.00,0.00){$t$}}
\put(200.00,20.00){\makebox(0.00,0.00){$\bullet$}}
\put(190.00,30.00){\line(0,-1){10.00}}
\put(190.00,20.00){\line(1,-2){10.00}}
\put(190.00,10.00){\line(-1,-3){3.33333}}
\put(210.00,10.00){\line(0,-1){10.00}}
\put(200.00,20.00){\line(1,-1){10.00}}
\put(200.00,20.00){\line(-1,-1){10.00}}
\put(200.00,30.00){\line(0,-1){10.00}}
\put(220.00,20.00){\line(0,-1){20.00}}
\put(210.00,30.00){\line(1,-1){10.00}}
\put(200.00,70.00){\line(0,-1){20.00}}
\put(200.00,50.00){\line(-1,-1){20.00}}
\put(220.00,30.00){\line(-1,1){20.00}}
\put(180.00,30.00){\line(1,0){40.00}}
%
\put(100.00,20.00){\makebox(0.00,0.00){$\bullet$}}
\put(120.00,40.00){\makebox(0.00,0.00){$t$}}
\put(90.00,20.00){\makebox(0.00,0.00){$\chi$}}
\put(130.00,30.00){\line(0,-1){30.00}}
\put(120.00,30.00){\line(0,-1){30.00}}
\put(107.00,10.00){\line(0,-1){10.00}}
\put(100.00,20.00){\line(2,-3){6.7}}
\put(100.00,20.00){\line(-2,-3){6.7}}
\put(93.00,10.00){\line(0,-1){10.00}}
\put(110.00,30.00){\line(-1,-1){10.00}}
\put(120.00,50.00){\line(0,1){20.00}}
\put(120.00,50.00){\line(-1,-1){20.00}}
\put(140.00,30.00){\line(-1,1){20.00}}
\put(100.00,30.00){\line(1,0){40.00}}
%
\put(20.00,70.00){\line(0,-1){20.00}}
\put(20.00,50.00){\makebox(0.00,0.00){$\bullet$}}

\put(10.00,60.00){\makebox(0.00,0.00){$\chi$}}
\put(6.8,30.00){\line(0,-1){30.00}}
\put(20.00,50.00){\line(-2,-3){13}}
\put(20.00,50.00){\line(2,-3){13}}
\put(-7,0){
\put(40.00,20.00){\makebox(0.00,0.00){$t$}}
\put(50.00,10.00){\line(0,-1){10.00}}
\put(40.00,10.00){\line(0,-1){10.00}}
\put(30.00,10.00){\line(0,-1){10.00}}
\put(40.00,30.00){\line(-1,-1){20.00}}
\put(60.00,10.00){\line(-1,1){20.00}}
\put(20.00,10.00){\line(1,0){40.00}}
}
\end{picture}}
\end{center}
\caption{\label{jj}%
Equation~(\protect\ref{s_ni_jsem_se_natrapil}) for $n=3$.}
\end{figure}

\begin{corollary}
The map $\ccc$ of~(\ref{dnes_jsem_mel_prednasku_ve_Filadelfii}) 
induces a canonical map (denoted again $\ccc$)
\begin{equation}
\label{kralicek}
\ccc : \calZ_\calP \to H^*(\calB_\calP,\delta_\calP),
\end{equation}
therefore $\HAP*$ is a natural $\calZ_\calP$-algebra.
\end{corollary}

From reasons which become clear later we call operations induced by
elements of $\calZ_\calP$ the {\em cup products\/}.  The following
proposition in which $\Lie$ is the operad for Lie algebras
(see~A.\ref{posloucham_Boroveho}) shows that the operad $\calZ_\calP$
is {\em always\/} nontrivial (provided $\calP \not = {\bf 1}$) while the
map~(\ref{kralicek}) is {\em never\/} monic.

\begin{proposition}
\label{lisasek}
The operad $\calZ_\calP$ contains the canonical element $\chi$.  There
exists a unique map $L: \ss \Lie \to \calZ_\calP$ that sends the
generator $\ss \lambda \in \ss \Lie(2)$ into $\chi \in
\calZ_\calP(2)$.  All elements in the image of this map are
$\delta_\calP$-cohomologous to zero in $\calB_\calP$.
\end{proposition}

\noindent
{\bf Proof.}
Recall~\cite[Corollary~2.2.9(b)]{ginzburg-kapranov:DMJ94} that, for
each quadratic operad $\calP$, there exits a morphism of operads $\Lie
\to \calP \ot \calP^!$ that takes the generator $\lambda \in \Lie(2)$
into the identity operator in $\calP(2) \ot \calP(2)^\# \cong \calP(2)
\ot \calP^!(2)$. Let $L : \ss{} \Lie \to \ss{}(\calP \ot \calP^!)$ be
the suspension of this morphism.  Let us prove, using
Theorem~\ref{poslouchal_jsem_Kryla}, that $\chi \in \calZ_\calP(2)$.
Equation~(\ref{s_ni_jsem_se_natrapil}) for $t = \chi$ reads
\[
\chi \circ_2 \chi + \chi \circ_1 \chi + (\chi \circ_2 \chi)(12) =0,
\]
which can be written, due to the symmetry~(\ref{porad_to_knika}) of
$\chi$, as the Jacobi identity for a degree $1$ ``multiplication''
$\chi$:
\begin{equation}
\label{zase_jsem_se_opil}
\chi \circ_1 \chi +  (\chi \circ_1 \chi)(123) + (\chi \circ_1 \chi)(132) =0,
\end{equation}
or, pictorially,
\[
{
\unitlength=1.00000pt
\thicklines
\begin{picture}(170.00,55.00)(-160.00,0)
\put(-110.00,30.00){\makebox(0.00,0.00){$+$}}
\put(-50.00,30.00){\makebox(0.00,0.00){$+$}}
\put(15.00,30.00){\makebox(0.00,0.00){$=\ 0$}}
\put(-120.00,0.00){\makebox(0.00,0.00){$3$}}
\put(-80.00,0.00){\makebox(0.00,0.00){$1$}}
\put(-40.00,0.00){\makebox(0.00,0.00){$2$}}
\put(-160.00,0.00){\makebox(0.00,0.00){$1$}}
\put(-60.00,0.00){\makebox(0.00,0.00){$2$}}
\put(-20.00,0.00){\makebox(0.00,0.00){$3$}}
\put(-140.00,0.00){\makebox(0.00,0.00){$2$}}
\put(-100.00,0.00){\makebox(0.00,0.00){$3$}}
\put(-0.00,0.00){\makebox(0.00,0.00){$1$}}
\multiput(-35.00,20.00)(-60,0){3}{\makebox(0.00,0.00)[br]{$\chi$}}
\multiput(-25.00,30.00)(-60,0){3}{\makebox(0.00,0.00)[br]{$\chi$}}
\put(-30.00,20.00){\makebox(0.00,0.00){$\bullet$}}
\put(-90.00,20.00){\makebox(0.00,0.00){$\bullet$}}
\put(-150.00,20.00){\makebox(0.00,0.00){$\bullet$}}
\put(-140.00,30.00){\makebox(0.00,0.00){$\bullet$}}
\put(-80.00,30.00){\makebox(0.00,0.00){$\bullet$}}
\put(-20.00,30.00){\makebox(0.00,0.00){$\bullet$}}
\put(-150.00,20.00){\line(0,1){0.00}}
\put(-140.00,10.00){\line(-1,1){10.00}}
\put(-140.00,50.00){\line(0,-1){20.00}}
\put(-140.00,30.00){\line(-1,-1){20.00}}
\put(-120.00,10.00){\line(-1,1){20.00}}
\put(-80.00,50.00){\line(0,-1){20.00}}
\put(-90.00,20.00){\line(1,-1){10.00}}
\put(-80.00,30.00){\line(-1,-1){20.00}}
\put(-60.00,10.00){\line(-1,1){20.00}}
\put(-20.00,10.00){\line(-1,1){10.00}}
\put(-20.00,30.00){\line(-1,-1){20.00}}
\put(-0.00,10.00){\line(0,1){0.00}}
\put(-20.00,30.00){\line(1,-1){20.00}}
\put(-20.00,50.00){\line(0,-1){20.00}}
\end{picture}}
\]

\noindent 
But~(\ref{zase_jsem_se_opil}) is satisfied, because $\chi = L(\ss{}
\lambda)$ by definition, and $\ss{} \lambda \in \ss{}\Lie (2)$
satisfies the same condition in $\ss \Lie$. The inclusion $\Im(L)
\subset \calZ_\calP$ follows from the fact that $\Im(L)$ is generated
by $\chi$ and that $\calZ_\calP$ is a suboperad of $\ss{}(\calP \ot
\calP^!)$.

Let us prove that all elements in the image of $L$ are
$\delta_\calP$-cohomologous to zero. Let $\ell \in \preLie(2)$ be the
generator of the quadratic operad $\preLie$ for pre-Lie algebras and
let $\circ := {\it pre}I(\ell) \in \calB^0_\calP(2)$, where ${\it
preI} : \preLie \to \calB_\calP$ is the map considered
in~(\ref{nasadil_jsem_retez}) at the end of
Section~\ref{jsem_otrok}. It is easy to verify that then
$\chi = \delta_\calP (\circ)$. This finishes the proof of
Proposition~\ref{lisasek}, because $\Im(L)$ is generated by $\chi$.%
\qed

Suppose that $\calP$ is the symmetrization of a non-$\Sigma$ operad
$\ucalP$.  Given $t \in \ss(\calP \ot \calP^!)(n)$ as in
Theorem~\ref{poslouchal_jsem_Kryla}, $\ccc(t) \in \calB_\ucalP(n)$ if
and only if $t$ belongs to the $\Sigma_n$-closure of $\ss(\ucalP \ot
\ucalP^!)(n)$ in $\ss(\calP \ot \calP^!)(n)$, that is, if $t =
\underline{t}\sigma$ for some $\underline{t} \in \ss(\ucalP \ot
\ucalP^!)(n)$ and $\sigma \in \Sigma_n$. In the following non-$\Sigma$
version of Theorem~\ref{poslouchal_jsem_Kryla}, $\uchi \in \ss
(\ucalP \ot \ucalP^!)(2)$ is the non-$\Sigma$ canonical element
introduced in Section~\ref{s2}.

\begin{figure}[t]
\begin{center}
{
\unitlength=1.000000pt
\begin{picture}(400.00,60.00)(0.00,10.00)
\thicklines
\put(20.00,70.00){\line(0,-1){20.00}}
\put(20.00,50.00){\makebox(0.00,0.00){$\bullet$}}
\put(10.00,60.00){\makebox(0.00,0.00){$\underline{\chi}$}}
\put(6.8,30.00){\line(0,-1){30.00}}
\put(20.00,50.00){\line(-2,-3){13}}
\put(20.00,50.00){\line(2,-3){13}}
\put(-7,0){
\put(40.00,20.00){\makebox(0.00,0.00){$\underline{t}$}}
\put(50.00,10.00){\line(0,-1){10.00}}
\put(40.00,10.00){\line(0,-1){10.00}}
\put(30.00,10.00){\line(0,-1){10.00}}
\put(40.00,30.00){\line(-1,-1){20.00}}
\put(60.00,10.00){\line(-1,1){20.00}}
\put(20.00,10.00){\line(1,0){40.00}}
}
\put(280.00,40.00){\makebox(0.00,0.00){$\underline{t}$}}
\put(310.00,20.00){\makebox(0.00,0.00){$\underline{\chi}$}}
\put(300.00,20.00){\makebox(0.00,0.00){$\bullet$}}
\put(330.00,40.00){\makebox(0.00,0.00){$=$}}
\put(290.00,10.00){\line(0,-1){10.00}}
\put(300.00,20.00){\line(-1,-1){10.00}}
\put(310.00,10.00){\line(0,-1){10.00}}
\put(300.00,20.00){\line(1,-1){10.00}}
\put(290.00,30.00){\line(1,-1){10.00}}
\put(280.00,30.00){\line(0,-1){30.00}}
\put(270.00,30.00){\line(0,-1){30.00}}
\put(280.00,50.00){\line(0,1){20.00}}
\put(280.00,50.00){\line(-1,-1){20.00}}
\put(300.00,30.00){\line(-1,1){20.00}}
\put(260.00,30.00){\line(1,0){40.00}}
\put(240.00,40.00){\makebox(0.00,0.00){$=$}}
\put(160.00,40.00){\makebox(0.00,0.00){$=$}}
\put(70.00,40.00){\makebox(0.00,0.00){$=$}}
\put(190.00,20.00){\makebox(0.00,0.00){$\underline{\chi}$}}
\put(90.00,20.00){\makebox(0.00,0.00){$\underline{\chi}$}}
\put(200.00,40.00){\makebox(0.00,0.00){$\underline{t}$}}
\put(120.00,40.00){\makebox(0.00,0.00){$\underline{t}$}}
\put(200.00,20.00){\makebox(0.00,0.00){$\bullet$}}
\put(100.00,20.00){\makebox(0.00,0.00){$\bullet$}}
\put(210.00,10.00){\line(0,-1){10.00}}
\put(200.00,20.00){\line(1,-1){10.00}}
\put(190.00,10.00){\line(0,-1){10.00}}
\put(200.00,20.00){\line(-1,-1){10.00}}
\put(200.00,30.00){\line(0,-1){10.00}}
\put(220.00,20.00){\line(0,-1){20.00}}
\put(210.00,30.00){\line(1,-1){10.00}}
\put(180.00,20.00){\line(0,-1){20.00}}
\put(190.00,30.00){\line(-1,-1){10.00}}
\put(200.00,70.00){\line(0,-1){20.00}}
\put(200.00,50.00){\line(-1,-1){20.00}}
\put(220.00,30.00){\line(-1,1){20.00}}
\put(180.00,30.00){\line(1,0){40.00}}
\put(130.00,30.00){\line(0,-1){30.00}}
\put(120.00,30.00){\line(0,-1){30.00}}
\put(110.00,10.00){\line(0,-1){10.00}}
\put(100.00,20.00){\line(1,-1){10.00}}
\put(90.00,10.00){\line(0,-1){10.00}}
\put(110.00,30.00){\line(-1,-1){20.00}}
\put(120.00,50.00){\line(0,1){20.00}}
\put(120.00,50.00){\line(-1,-1){20.00}}
\put(140.00,30.00){\line(-1,1){20.00}}
\put(100.00,30.00){\line(1,0){40.00}}
\put(400,0){
\put(-20.00,70.00){\line(0,-1){20.00}}
\put(-20.00,50.00){\makebox(0.00,0.00){$\bullet$}}
\put(-10.00,60.00){\makebox(0.00,0.00){$\underline{\chi}$}}
\put(-6.8,30.00){\line(0,-1){30.00}}
\put(-20.00,50.00){\line(2,-3){13}}
\put(-20.00,50.00){\line(-2,-3){13}}
\put(7,0){
\put(-40.00,20.00){\makebox(0.00,0.00){$\underline{t}$}}
\put(-50.00,10.00){\line(0,-1){10.00}}
\put(-40.00,10.00){\line(0,-1){10.00}}
\put(-30.00,10.00){\line(0,-1){10.00}}
\put(-40.00,30.00){\line(1,-1){20.00}}
\put(-60.00,10.00){\line(1,1){20.00}}
\put(-20.00,10.00){\line(-1,0){40.00}}
}
}
\end{picture}}
\end{center}
\caption{\label{iii}
Equation~(\protect\ref{kocicka}) for $n=3$.}
\end{figure}
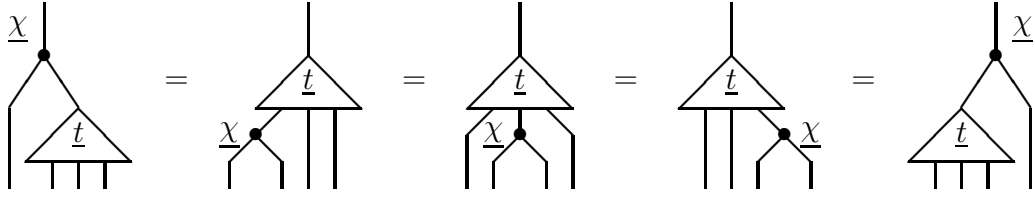

\begin{theorem}
\label{tuzka}
An element $\ut \in \ss(\ucalP \ot \ucalP^!)(n)\subset \ss(\calP \ot
\calP^!)(n)$ belongs to
$\calZ_\calP(n)$ if and only if
\begin{equation}
\label{kocicka}
\uchi \circ_2 \ut = \ut \circ_1 \uchi =  \ut \circ_2 \uchi
= \cdots =  \ut \circ_n \uchi= \uchi \circ_1 \ut,
\end{equation}
see Figure~\ref{iii}.
\end{theorem}

{\bf Proof} of Theorem~\ref{tuzka} 
is a straightforward verification. The {\bf proof}
of the following proposition is similar to that of
Proposition~\ref{lisasek}.

\begin{proposition}
\label{posledni_den_v_Phily}
Let $\calP$ be the symmetrization of a non-$\Sigma$ operad
$\ucalP$. Then $\uchi \in \calZ_\calP$ and there exists a unique map
$A : \ss\As \to \calZ_\calP$ defined by $A(\ss \mu) := \uchi$, where
$\ss\mu \in \ss\As(2)$ is the suspension of the generator $\mu$
(see~A.\ref{posloucham_Boroveho}). Moreover, the diagram
\[
\Dtriangle {\ss\Lie}{\calZ_\calP}{\ss\As}{L}{}{A}
\]
where $L$ is as in Proposition~\ref{lisasek},
with the vertical map given by the anti-commutator of the
associative product, commutes.
\end{proposition}

\begin{example}
{\em -- Hochschild cohomology.
Let $\calP = \Ass{}$ be the operad for associative algebras. 
Then $\ss(\calP \ot \calP^!) = \ss(\As \ot \As)$
and a simple calculation reveals that the map $A$ of
Proposition~\ref{posledni_den_v_Phily} is the suspended diagonal $\ss\Delta:
\ss\As \to \ss(\As \ot \As)$ and that $\calZ_\Ass = \Im(A)$. Therefore 
\[
\calZ_\As \cong \ss \As.
\] 
The generator $\ss\mu\in \ss \As(2)$ is mapped to the ``classical'' cup product
$f,g \mapsto f\cup g$ of Hochschild cochains~\cite{gerstenhaber:AM63}, 
and the generator $\ss\lambda \in \ss\Lie(2)$ to the anti-commutator
of this cup product:
\[
f,g \mapsto f\cup g +(-1)^{|f||g|} g \cup f
\]
which cohomologous to zero, because the cup product of Hochschild
cochains is homotopy commutative~\cite[Theorem~3]{gerstenhaber:AM63}.
}  
\end{example}

\begin{example}
\label{oslicek1} 
{\em -- Chevalley-Eilenberg cohomology. If
$\calP = \Lie$ is the operad for Lie algebras, 
then $\ss(\calP \ot \calP^!) = \ss(\Lie \ot \Com) \cong
\ss\Lie$
and we see immediately that
\begin{equation}
\label{hhh}
\calZ_\Lie \cong \ss\Lie = \Im(L).
\end{equation}
The generator $\ss\lambda \in \Lie(2)$ is mapped to the product $f,g
\mapsto \{f,g\}$, which is cohomologous to zero,
see~\cite[Exercise~7]{lm:sb}. 
}  
\end{example}

\begin{example}
\label{oslicek2}
{\em -- Harrison cohomology. Here
$\calP = \Com$ is the operad for commutative associative algebras 
and  $\ss(\calP \ot \calP^!) = \ss(\Com \ot \Lie)  \cong \ss\Lie$,
therefore, as in Example~\ref{oslicek1},
\begin{equation}
\label{1hhh}
\calZ_\Com \cong \ss \Lie= \Im(L).
\end{equation}
Equations~(\ref{hhh}) and~(\ref{1hhh}) illustrate the obvious
self-duality of the space of cup products:
\[
\calZ_{\calP^!} \cong \calZ_{\calP},
\]
compare Conjecture~\ref{zase}.
}  
\end{example}

\begin{example}
{\em
If $\calD = \As * \As$ is as in 
Example~\ref{ve_stredu_poletim_do_Minneapolis}, then
\[
\calZ_\calD = \ss \As \vee \ss \As.
\]
Let us describe products corresponding to the
generators of the $4$-dimensional vector space
\[
(\ss \As \vee \ss \As)(2) =
\ss\As(2) \oplus \ss\As(2) 
\cong \hskip .2em \uparrow \bfk[\Sigma_2] \oplus  \uparrow
\bfk[\Sigma_2].
\]

Recall that $\CAD*$ is the total complex of the meager bicomplex
in Figure~\ref{jsem_nachlazeny}.
Let $\cup^1$ (resp.~$\cup^2$) be the cup product
in the horizontal (resp.~vertical) subcomplex in
Figure~\ref{jsem_nachlazeny}.  
Let $\pi_1$ (resp.~$\pi_2$) be the projection of $\CAD*$ onto the horizontal
(vertical) subcomplex. Likewise, let $\iota_1$
(resp.~$\iota_2$) be the inclusion. Although neither
$\pi_i$, $\iota_i$ nor $\cup^i$ are chain maps  ($i =
1,2$), the compositions
\[
f \cup_1 g : = \iota_1(\pi_1 f \cup^1 \pi_1 g)\
\mbox{ and } 
f \cup_2 g : = \iota_2(\pi_2 f \cup^2 \pi_2 g)
\]
are chain operations. The generators of $\ss\As(2) \oplus \ss\As(2)$
then correspond to the four operations 
\[
f,g \mapsto f\cup_1 g,\  f,g \mapsto f\cup_2 g,\
f,g \mapsto g\cup_1 f\ \mbox { and } f,g \mapsto g\cup_2 f.
\] 
The combination
\[
(f\cup_1 g + f\cup_2 g) + (-1)^{|f||g|} (g\cup_1 f + g\cup_2 f)
\]
is cohomologous to zero and the image $T(\calZ_\calD(2))$ of
$\calZ_\calD$ in $H^1(\calB_\calD(2))$ is easily seen to be
$3$-dimensional.  }
\end{example}

\section{Operad $\calB_\Ass$ and the Deligne conjecture}
\label{pisi_na_novem_notebooku}
\label{s6}

In this section we recall some results related to $\calB_\Ass$ and the
Deligne conjecture. Let us make a necessary comment about our degree
convention. We use the grading such that the intrinsic bracket
of Section~\ref{s4} has degree $0$ in $\calB^*_\calP(2)$, while the
$n$-fold cup products of Section~\ref{s5} are of degree $n-1$ in
$\calB^*_\calP(n)$. In the literature related to the Deligne
conjecture, the convention under which the intrinsic bracket has degree
$1$ and the $n$-fold cup products are of degree
$0$ is often used. These two conventions are
tied by the following regrading operator:
\[
\Reg(\calB^*_\calP (n)) := \calB^{n-1-*}_\calP(n).
\]

In what follows we identify operads that differ only by the above
regrading. In particular, the operad $\Gerst$ for Gerstenhaber algebras
becomes identified with the operad $\Braid$ for braid algebras (also
called 1-algebras), see~A.\ref{ma_prejit_fronta_v_sobotu}. 

Let us recall that a topological operad $\calA$ is an {\em
$E_2$-operad\/} if it has the homotopy type of the little discs operad
${\disc}_2$~\cite{may:1972}. According to the Formality
Theorem~\cite{tamarkin:formality}, the operad $S_*(\calA)$ of singular
chains on such an operad has the homotopy type of the operad 
$\Braid$ for braid algebras.

(i) {\em D.~Tamarkin\/} and {\em B.~Tsygan\/} studied
in~\cite[Section~3]{tamarkin-tsygan:LMP01} a certain operad $F =
\{F(n)\}_{n \geq 1}$ of natural operations on the cosimplicial
Hochschild complex $C^\bullet(X,X)$ of a topological unital
monoid~$X$. The $n$-th space of this cosimplicial set is the space
${\it Cont\/}(X^{\times n},X)$ of continuous maps from the $n$-th
cartesian power of $X$ to $X$. For each $n \geq 1$, $F(n)$ is
a functor $(\Delta^{\it op})^n \times \Delta \to
{\sf Sets}$. They then considered a topological operad $E =
\{E(n)\}_{n \geq 1}$ whose pieces are the topological realizations of
these functors and claimed that $E$ is an $E_2$-operad. 

It is not difficult to see that the operad ${\it CN\/}_*(F)$ of
normalized chains of $F$ coincides with the operad
$\{\calB_\uAss(n)\}_{n \geq 1}$ (our $\calB_\uAss$ without
constants). Since $(\calB^*_\uAss(0),\delta_\uAss)$ is acyclic (see
Example~\ref{ctvrek_a_zase_je_lat}), we could conclude that
$\calB_\uAss$ has the homotopy type of $\Gerst$, but we must bear in
mind that the arguments in~\cite{tamarkin-tsygan:LMP01} were merely
sketched.

(ii) {\em J.E.~McClure\/} and {\em J.H.~Smith\/} considered
in~\cite{mcclure-smith:JAMS03} a dg-suboperad $\calS_2$ of their
``sequence'' operad $\calS$ and proved that $\calS_2$ naturally acts
on the Hochschild cochain complex of an associative algebra. In our
terminology this means that they constructed a canonical map $\calS_2
\to \calB_\Ass{}$. They then verified the Deligne conjecture by showing,
using a result of~\cite{berger:CM202}, that $\calS_2$ has the homotopy
type of the singular chain complex $S_* (\disc_2)$ of the little discs
operad.  Their proof is a very reliable one.

(iii) {\em M.~Kontsevich\/} and {\em
Y.~Soibelman\/}~\cite{kontsevich-soibelman} introduced a ``minimal
operad'' $M$ naturally acting on the Hochschild cochain complex of an
$A_\infty$-algebra. In our terminology, $M$ was a suboperad of
$\calB_{sh\Ass}$ generated by braces and cup products.  They then
argued that $M$ has the homotopy type of the operad of suitably
defined piecewise algebraic chains on the operad
$FM_2$ of the Fulton-MacPherson compactification of the configuration space of
points in ${\mathbb R}^2$. Since $FM_2$ is,
by~\cite[Proposition~3.9]{salvatore}, an $E_2$-operad,
they concluded that $M$ has the homotopy type of $\Gerst$. 

(iv) {\em R.M.~Kaufmann\/} realized in~\cite{kaufmann}
that the cellular chains ${\it CC\/}_*(\Cact^1)$ on his operad-up-to-homotopy
$\Cact^1$ of spineless normalized cacti is a honest operad which
naturally acts on the Hochschild cochain complex,
via braces and cup products. By comparing $\Cact^1$ to the operad
$\Cact$ of spineless (non-normalized) cacti, he concluded that
$CC_*(\Cact^1)$ is a model for chains on the little discs operad
$\disc_2$.

All the proofs of the Deligne conjecture mentioned above use
some special features of associative algebras and
$E_2$-operads, such as the cosimplicial structure of the Hochschild
cochain complex, Fiedorowicz' detection principle, or a relation to the
Fulton-MacPherson and cacti operads. None of these features are
available for a general operad $\calP$, we therefore think that the
analysis of the homotopy type of $\calB_\calP$ for a general Koszul
quadratic $\calP$ is substantially more difficult than the analysis of
$\calB_\As$. 

Let us mention that there are other approaches to the Deligne
conjecture, as D.~Tamarkin's proofs that use the Etingof-Kazhdan
quantization~\cite{hinich:FM03,tamarkin:deligne}, or those
based on a suitable filtration of the Fulton-MacPherson
compactification $FM_2$, see E.~Getzler and
J.D.S.~Jones~\cite{getzler-jones:preprint} or
A.A.~Voronov~\cite{voronov:99}.

\section{Natural operations}
\label{natural}

Let use recall the following definitions which can be found for
example in~\cite[Section~II.1.5]{markl-shnider-stasheff:book}.  By a
{\em tree\/} we mean a connected graph $T$ without loops. A {\em
valence\/} of a vertex $v$ of $T$ is the number of edges adjacent to
$v$. A {\em leg\/} or {\em leaf\/} of $T$ is an edge adjacent to a
vertex of valence one, other edges of $T$ are {\em interior\/}.  We in
fact discard vertices of valence one at the endpoints of the legs,
therefore the legs become ``half-edges'' having only one vertex. By a
{\em rooted\/} or {\em directed\/} tree we mean a tree with a
distinguished {\em output\/} leg called the {\em root\/}.  The
remaining legs are called the {\em input legs\/} of the tree. A tree
with $a$ input legs labelled by elements of the set $\{1,2,\ldots,a\}$
is called an {\em $a$-tree\/}.  A rooted tree is automatically {\em
oriented\/}, each edge pointing towards the root.  The edges pointing
towards a given vertex $v$ are called the {\em input edges\/} of $v$,
the number of these input edges is then the {\em arity\/} of $v$
denoted ${\it ar\/}(v)$. Vertices of arity one are called unary,
vertices of arity two binary, vertices of arity three ternary,~etc.

\noindent 
{\bf Notation.}
{\it\
Let $n,m$ and $m_1,\ldots,m_n$ be non-negative integers. 
In the rest of this section, $i$ will always denote an integer between
$0$ and $n$,  $a:= m+1$ and $a_i := m_i + 1$. We will also assume the
notation introduced in~A.\ref{Ferda}.\/}

An {\em $n$-linear natural operation\/} 
\[
U : \CAP{m_1} \ot \cdots \ot \CAP{m_n} \to \CAP m
\]
is given by the following data.

(i) A rooted $a$-tree $T$ with $n$ white vertices
$w_1,\ldots,w_n$ of  arities $a_1,\ldots,a_n$,
and $k$ at least binary black vertices, $k \geq 0$. 

(ii) A linear order on the set of input edges of each white vertex of $T$.

(iii) A decoration of black vertices of $T$ by elements of $\calP$.

(iv) A linear map $\Phi : \ss \calP^!(a_1) \ot \cdots \ot
\ss \calP^!(a_n)  \to \ss \calP^!(a)$.

Given the above data and $f_i \in \CAP{m_i}$, the value
$U(f_1,\ldots,f_n) \in \CAP m$ is defined as follows. Let us decompose
\[
f_i = \sum_{\kappa_i} \phi_i^{\kappa_i} \ot q^i_{\kappa_i} \in [\Lin
  (A^{\otimes a_i},A) \ot \ss \calP^!(a_i)]^{\Sigma_{a_i}} \cong
\CAP{m_i},
\]
where $\phi_i^{\kappa_i} \in \Lin(A^{\ot a_i},A)$, $q^i_{\kappa_i} \in
\ss\calP^!(a_i)$ and $\kappa_i$ is a summation index.  Since the
inputs of white vertices are linearly ordered, each
$\phi_i^{\kappa_i}$ determines a decoration of the white vertex $w_i$
by an element of $\Lin( A ^{\ot a_i}, A) = \End_A(a_i)$. Recall that
$A$ is a $\calP$-algebra with the structure homomorphism $\alpha :
\calP \to \End_A$. Applying $\alpha$ to the decorations of the black
vertices we decorate also black vertices with elements of $\End_A$. So
$T$ is now a tree with all vertices decorated by $\End_A$.  The
composition in the operad $\End_V$ along
$T$~\cite{ginzburg-kapranov:DMJ94} determines, for each
$k_1,\ldots,k_n$, the element
\[
T(\phi_1^{\kappa_1},\ldots,\phi_n^{\kappa_n}) \in \Lin(A^{\ot a}, A).
\]
Let 
\[
\widetilde U(f_1,\ldots,f_n) := \sum_{\kappa_1,\ldots,\kappa_n} 
T(\phi_1^{\kappa_1},\ldots,\phi_n^{\kappa_n}) \ot
\Phi(q^1_{\kappa_1},\ldots,q^n_{\kappa_n}) \in \Lin
 (A^{\ot a},A) \ot \ss\calP^!(a) \cong \tildeCAP m.
\]
Finally, let $U(f_1,\ldots,f_n) := {\it Aver\/}(\widetilde
U(f_1,\ldots,f_n)) \in \CAP m$.
It follows from an elementary combinatorics of trees that
\[
\deg ( U) = {\it ar\/}(b_1) + \cdots + {\it ar\/}(b_k) - k,
\]
therefore $\deg(U)$ is {\em always non-negative\/} and
$\deg(U) = 0$ if and only if $T$ has no black vertex.

\begin{definition}
\label{sam_s_detickami}
Let $\calB_\calP := \{\calB_\calP(n)\}_{n \geq 0}$ be the operad
spanned by all natural operations $U = U_{(T,\Phi)}$ in the above sense.
Since the differential $d_\calP$ of $\CAP*$ is itself a natural
operation living in $\calB_\calP^1(1)$, it induces a differential
$\delta_\calP$ on $\calB_\calP$ by the standard formula
\[
\delta_\calP(U)(f_1,\ldots,f_n) := d_\calP U(f_1,\ldots,f_n) -
(-1)^{|U|} \sum_{1 \leq i \leq n} 
(-1)^{|f_1| + \cdots + |f_{i-1}|}
U(f_1,\ldots,d_\calP f_i,\ldots, f_n),
\]
making $\calB_\calP = (\calB^*_\calP,\delta_\calP)$ a dg-operad.
\end{definition}

Heuristically, the value $U_{(T,\Phi)}(f_1,\ldots,f_n)$ is given by
inserting $f_i$ at the vertex $w_i$ of $T$, $1\leq i \leq n$, and then
performing the composition along $\Phi$. The operadic composition of
$\calB_\calP$ is the vertex insertion similar to that
of~\cite{chapoton-livernet:pre-lie} and the symmetric group permutes
the labels of white vertices. In the following definition we introduce
a non-$\Sigma$ version of $\calB_\calP$.

\begin{definition}
\label{letel_jsem_do_Nymburka!!}
Suppose $\calP$ is the symmetrization of a non-$\Sigma$ operad
$\ucalP$. Let $\calB_\ucalP$ be the dg-suboperad of $\calB_\calP$
spanned by natural operations $U_{(T,\Phi)}$ as in
Definition~\ref{sam_s_detickami} such that the tree $T$ is
planar, with black vertices decorated by elements of \hskip .1em $\ucalP$,
and the map $\Phi$ such that 
\[
\Phi(\ss \ucalP^!(a_1)
\ot \cdots \ot \ss \ucalP^!(a_n)) \subset \ss \ucalP^!(a).
\]
\end{definition}

\begin{example}
{\rm
--\label{Kocicka_zrcatko}
Constants. Let us see what happens if $T$ is the
$a$-corolla with one black vertex decorated by $p \in \calP(a)$ and no
white vertices as in Figure~\ref{rozbil_jsem_notebook}. 
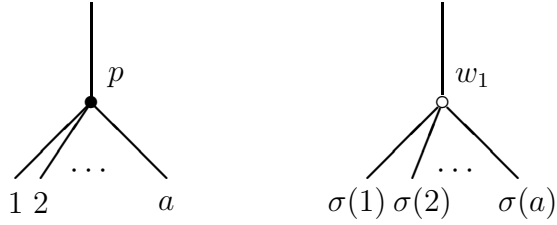
\begin{figure}[t]
\begin{center}
\unitlength=.95pt
\begin{picture}(60.00,80.00)(20.00,0.00)
\thicklines
\put(-70,0){
\put(0.00,0.00){\makebox(0.00,0.00){$1$}}
\put(10.00,0.00){\makebox(0.00,0.00){$2$}}
\put(60.00,0.00){\makebox(0.00,0.00){$a$}}
\put(30.00,10.00){\makebox(0.00,0.00)[b]{$\cdots$}}
\put(40.00,50.00){\makebox(0.00,0.00){$p$}}
\put(30.00,40.00){\makebox(0.00,0.00){$\bullet$}}
\put(60.00,10.00){\line(-1,1){30.00}}
\put(10.00,10.00){\line(2,3){20.00}}
\put(0.00,10.00){\line(1,1){30.00}}
\put(30.00,80.00){\line(0,-1){40.00}}
}
\put(70,0){
\put(42.00,50.00){\makebox(0.00,0.00){$w_1$}}
\put(-4.00,0.00){\makebox(0.00,0.00){$\sigma(1)$}}
\put(22.00,0.00){\makebox(0.00,0.00){$\sigma(2)$}}
\put(64.00,0.00){\makebox(0.00,0.00){$\sigma(a)$}}
\put(36.00,10.00){\makebox(0.00,0.00)[b]{$\cdots$}}
\put(30.00,40.00){\makebox(0.00,0.00){$\circ$}}
\put(60.00,10.00){\line(-1,1){28.00}}
\put(18.00,10.00){\line(2,5){11.00}}
\put(0.00,10.00){\line(1,1){28.00}}
\put(30.00,80.00){\line(0,-1){37.00}}
}
\end{picture}
\end{center}
\caption{\label{rozbil_jsem_notebook}
The tree defining a constant in $\calB^m_\calP(0)$ 
(left) and a unary operation in $\calB^0_\calP(1)$ (right), where $m =
a+1$ as always.}
\end{figure}
The map
$\Phi: \bfk \to \ss \calP^!(a)$ is given by specifying an element
$\varphi := \Phi(1) \in \ss\calP^!(a)$ and  $\widetilde U$ determined
by this $\Phi$ equals $\alpha
(p) \ot \varphi \in \tildeCAP m$. Since $\alpha$ is equivariant,
\[
{\it Aver\/}(\alpha(p) \ot \varphi) = (\alpha \ot \id)({\it Aver\/}(p \otimes
\varphi)) 
\]
therefore $U := {\it Aver\/}(\alpha(p) \ot \varphi) \in \CAP m$ 
is parametrized by an element in the image of the averaging map 
\[
{\it  Aver\/}: \calP(a) \ot \ss \calP^!(a) \to (\calP(a) \ot \ss
\calP^!(a))^{\Sigma_a},
\] 
in other words, 
\[
\calB_\calP^{m}(0) \cong \ss (\calP \ot \calP^!)(a)^{\Sigma_a},\ m \geq 0.
\]
It is equally easy to see that, for a quadratic Koszul non-$\Sigma$
operad $\ucalP$,
\[
\calB^m_\ucalP(0) \cong \ss (\ucalP \ot\ucalP^!)(a),\ m \geq 0.
\]
}
\end{example}

\begin{example}
-- \label{Deda_Mraz} {\rm\ Unary operations of degree $0$. Now $T$ is
  an $a$-corolla with one white planar vertex and no black vertices,
  with input legs labelled $\sigma(1),\ldots,\sigma(a)$, $\sigma \in
  \Sigma_a$, as shown in Figure~\ref{rozbil_jsem_notebook}, and $\Phi :
  \ss\calP^!(a) \to \ss\calP^!(a)$ is a linear map. If
\[
f = \sum_\kappa \phi^\kappa \ot q_\kappa \in [\Lin( \otexp Aa ,A) \ot
  \ss \calP^!(a)]^{\Sigma_a} \cong \CAP m,
\]
then 
$U(f) = {\it Aver\/}(\sum_\kappa\phi^\kappa \sigma^{-1} \ot \Phi(q_\kappa))$.
Since $f = \sum \phi^\kappa \ot q_\kappa$ is $\Sigma_a$-stable,
\[
U(f) =
\sum \phi^\kappa \ot {\it Aver\/}(\Phi_\kappa)(q_\kappa).
\]
Therefore $U(f)$ is given  by a $\Sigma_a$-equivariant map 
$\Psi := \desusp^{a-1}{\it Aver\/}(\Phi) : \calP^!(a) \to \calP^!(a)$,
thus
\begin{equation}
\label{ll}
\calB_\calP^0(1) \cong \Lin_\Sigma(\calP^!,\calP^!),
\end{equation}
the space of all collections $\{\psi_n: \calP^!(n) \to \calP^!(n)\}_{n
\geq 0}$ of equivariant maps. We leave as an exercise to verify that,
for a non-$\Sigma$ quadratic Koszul operad $\ucalP$,
\[
\calB_\ucalP^0(1) \cong \Lin(\ucalP^!,\ucalP^!).
\]
}
\end{example}

\begin{example}
-- {\rm 
Projections. 
Let $p_m \in \calB_\calP^0(1)$ be given, in identification~(\ref{ll}), by
$\Psi \in \Lin_\Sigma(\calP^!,\calP^!)$ defined as
\[
\Psi|_{\calP^!(a)} =
\cases{\id_{\calP^!(a)}}{for $a=m+1$ and}0{othewise.}
\] 
Clearly, $p_m$ is the projection $\CAP* \epi \CAP
m$. The system of all these
projections makes $\calB_\calP$ an $\BbbZ^{\geq 0}$-colored operad, where
$\BbbZ^{\geq 0}$ is the set of non-negative integers. Since 
these projections do not commute with $d_\calP$ (that is
$\delta_\calP(p_m) \not= 0$ for a generic $\calP$),
$(\calB_\calP,\delta_\calP)$ {\em is not\/} a {\em dg\/}
$\BbbZ^{\geq 0}$-colored operad.
}
\end{example}

\begin{example}
\label{fit}
-- {\rm Cup products.  In this example we explain how an element 
\[
t = p
  \ot \ss q \in \calP(n) \ot \ss\calP^!(n) \cong \ss(\calP \ot
  \calP^!)(n)
\] 
determines a natural operation in $\calB_\calP^{n-1}(n)$.  Let $T$ be
as in Figure~\ref{dnes_novy_notebook}, with the black vertex
decorated by $p \in \calP(n)$,
\begin{figure}[t]
\begin{center}
{
\thicklines
\unitlength=.8pt
\begin{picture}(240.00,100.00)(0.00,10.00)
\put(220.00,0.00){\makebox(0.00,0.00){$\underbrace{\rule{33pt}{0pt}}_{a_n}$}}
\put(70.00,0.00){\makebox(0.00,0.00){$\underbrace{\rule{33pt}{0pt}}_{a_2}$}}
\put(20.00,0.00){\makebox(0.00,0.00){$\underbrace{\rule{33pt}{0pt}}_{a_1}$}}
\put(130.00,30.00){\makebox(0.00,0.00){$\cdots$}}
\put(224.00,10.00){\makebox(0.00,0.00)[b]{\scriptsize $\cdots$}}
\put(74.00,10.00){\makebox(0.00,0.00)[b]{\scriptsize $\cdots$}}
\put(24.00,10.00){\makebox(0.00,0.00)[b]{\scriptsize $\cdots$}}
\put(130.00,100.00){\makebox(0.00,0.00){$p$}}
\put(220.00,30.00){\makebox(0.00,0.00){$\circ$}}
\put(70.00,30.00){\makebox(0.00,0.00){$\circ$}}
\put(20.00,30.00){\makebox(0.00,0.00){$\circ$}}
\put(10.00,35.00){\makebox(0.00,0.00)[rb]{$w_1$}}
\put(50.00,35.00){\makebox(0.00,0.00)[lb]{$w_2$}}
\put(230.00,35.00){\makebox(0.00,0.00)[lb]{$w_n$}}
\put(120.00,90.00){\makebox(0.00,0.00){$\bullet$}}
\put(120.00,120.00){\line(0,-1){30.00}}
\put(222,28){\line(1,-1){18}}
\put(218.8,27.5){\line(-1,-2){8.5}}
\put(218,28){\line(-1,-1){18}}
\put(-200,0){
\put(222,28){\line(1,-1){18}}
\put(218.8,27.5){\line(-1,-2){8.5}}
\put(218,28){\line(-1,-1){18}}
}
\put(-150,0){
\put(222,28){\line(1,-1){18}}
\put(218.8,27.5){\line(-1,-2){8.5}}
\put(218,28){\line(-1,-1){18}}
}
\put(120.00,90.00){\line(5,-3){97}}
\put(120.00,90.00){\line(-5,-6){48}}
\put(120.00,90.00){\line(-5,-3){97}}
\end{picture}}
\end{center}
\caption{\label{dnes_novy_notebook}
The tree defining the cup product.}
\end{figure}
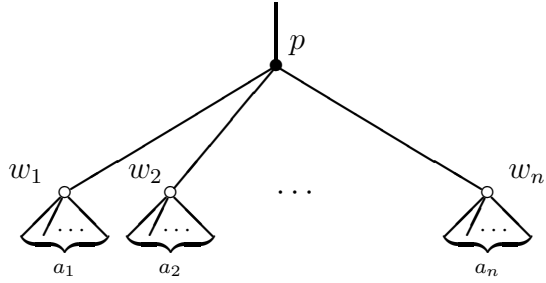
and let the linear map 
$\Phi : \ss\calP^!(a_1) \ot \cdots \ot
\ss\calP^!(a_n) \to \ss \calP^!(a)$ be given by the operadic composition in
$\ss{}\calP^!$:
\[
\Phi(\ss q_1,\ldots,\ss  q_n) :=
\ss q(\ss q_1,\ldots,\ss q_n),\ q_i \in \calP^!(a_i),\ 1 \leq i \leq n.
\]
It is more or less clear that the natural operation $U_{(T,\Phi)}$ 
determined by the above data
agrees with the cup product $\ccc(t)$ of
Theorem~\ref{zitra_prednaska}. We recommend as another exercise to verify
that also the intrinsic bracket described in~(\ref{lepidlo}) is
given by natural operations in the sense of this section.  
}
\end{example}

\begin{example}
{\rm
Let us describe all natural operations $\CAP 1 \ot \CAP 1 \to \CAP 2$
for some particular choices of $\calP$.

(i) Hochschild cohomology.  For $\calP = \As$, $\CAP 1 =
\Lin (\otexp A2,A)$, $\CAP 2 =  \Lin (\otexp A3,A)$, and the
only natural operations $\CAP 1 \ot \CAP 1 \to \CAP 2$ are linear
combinations of
\[
f,g \to (f \circ_1 g)\sigma,\ f,g \to (f \circ_2 g)\sigma,\ f,g \to 
(g \circ_1 f)\sigma,\ f,g \to (g \circ_2 f)\sigma,\ \sigma \in \Sigma_3,
\]
where $\circ_1$, $\circ_2$ are Gerstenhaber-type
products~\cite{gerstenhaber:AM63} given by
\begin{equation}
\label{pisu_v_susarne}
(u \circ_1 v) (a,b,c) := u(v(a,b),c)),\
(u \circ_2 v) (a,b,c) := u(a,v(b,c)),\ 
\end{equation}
for $u,v \in \CAP 2$, $a,b,c \in A$, and $\sigma \in \Sigma_3$
permutes the factors of $\otexp A3$ in the usual way. Operations
belonging to $\calB^{0}_\uAss(2)$ are linear combinations of the
operations~(\ref{pisu_v_susarne}) with $\sigma
= \id_3$, the unit of $\Sigma_3$.

(ii) Chevalley-Eilenberg cohomology.
If $\calP = \Lie$, then 
 $\CAP 1 =  \Lin (\land^2,A)$, $\CAP 2 =
\Lin (\land^3 A,A)$, where $\land^n A$ denotes the $n$-th
exterior power. The only natural operations  $\CAP 1
\ot \CAP 1 \to \CAP 2$ are linear combinations of
\[
f,g \to f \circ g\ \mbox { and } \ f,g \to g \circ f,
\]
where
\[
(u\circ v)(a,b,c) := u(v(a,b),c) +  u(v(b,c),a) + u(v(c,a),b)
\] 
for $u,v \in \CAP 2$ and $a,b,c \in A$.

(iii) Harrison cohomology.
If $\calP = \Com$, then 
\[
\CAP 1 =  \{u \in   \Lin (\otexp A2,A);\
u(a,b) - u(b,a) = 0\}
\]  
and $\CAP 2$ consists of all 
$w \in \Lin (\otexp A3,A)$ such that 
\[
w(a,b,c) - w(b,a,c) + w(b,c,a) =
w(a,b,c) - w(a,c,b) + w(c,a,b) = 0,
\]
for $a,b,c \in A$.
Natural operations $\CAP 1 \ot \CAP 1 \to \CAP 2$ 
are linear combinations of
\[
f,g \to f \circ g\ \mbox { and } \ f,g \to g \circ f,
\]
where
\[
u \circ v := u(v(a,b)c) - u(v(b,c),a),
\] 
for $u,v \in \CAP 2$ and $a,b,c \in A$.
}
\end{example}

\section*{Appendix: Notations, conventions and background material}

\begin{Appendix}
\label{jitka_neodpovida}
{\rm
In this note, an {\em operad\/} means an operad in the category of
differential graded (dg) vector spaces, that is, a sequence 
$\calP = \{\calP(n)\}_{n \geq 0}$ of right $\Sigma_n$-modules with structure
operations 
\[
\gamma : \calP(n) \ot \calP(k_1) \ot \cdots \ot \calP(k_n) 
\to
\calP(k_1 + \cdots + k_n),
\]
for $n \geq 1$ and $k_1,\ldots,k_n \geq 0$, and a unit map $\eta :
\bfk \to \calP(1)$ that satisfy the usual
axioms~\cite{may:1972,kriz-may}. Instead of $\gamma(p \ot p_1 \ot
\ldots \ot p_n)$ we will often write $\gamma(p, p_1, \ldots , p_n)$ or
$p(p_1,\ldots,p_n)$. Recall~\cite{markl:zebrulka} that operads can be
equivalently defined using the $\circ_i$-operations
\[
\circ_i : \calP(m) \ot \calP(n) \to \calP(m+n-1)
\] 
defined, for $m,n \geq 0$, $1 \leq i \leq m$, by
\[
p\circ_i q := \gamma(p \ot e^{\ot (i-1)}\ot q \ot \otexp e{m-i}), 
\]
where $e : = \eta(1)$.

If we remove from the above definition all references to the symmetric group
actions, we get a definition of a {\em non-$\Sigma$ operad\/}. Each
non-$\Sigma$ operad $\ucalP$ generates a unique (usual) operad $\calP$
such that $\calP(n) \cong \ucalP(n) \ot \bfk[\Sigma_n]$, $n \geq 0$.
}
\end{Appendix}

\begin{Appendix}
\label{posloucham_Boroveho}
{\rm For each set of operations $E$, there exists the {\em free
operad\/} $\Gamma(E)$ generated by
$E$~\cite[Proposition~II.1.92]{markl-shnider-stasheff:book}.  Let
$\mu$ denote a bilinear operation placed in degree $0$. The operad
$\Ass{}$ for {\em associative algebras\/} is the quotient
\[
\Ass{} := \Gamma(\mu)/(\mu \circ_1 \mu - \mu \circ_2 \mu),
\]
where $(\mu \circ_1 \mu - \mu \circ_2 \mu)$ denotes the operadic ideal
generated by the associativity axiom for $\mu$.

If $\lambda$ is a skew-symmetric bilinear operation,
then the operad for {\em Lie algebras\/} is the quotient
\[
\Lie := \Gamma(\lambda)/({\it Jacobi}(\lambda)),
\]
where 
\[
{\it Jacobi}(\lambda) := \sum_{\sigma\in C_3} (\lambda \circ_1 \lambda) \sigma
\]
with the summation taken over the order $3$ cyclic subgroup $C_3$ of
$\Sigma_3$, denotes the Jacobi identity for $\lambda$.

Finally, for an arbitrary differential graded vector space $V$, there
is the {\em endomorphism operad\/} $\End_V = \{Lin(V^{\ot n},V)\}_{n
\geq 0}$, with structure operations given as the usual composition of
multilinear maps. A {\em $\calP$-algebra\/} is then a homomorphism
$\alpha : \calP \to \End_V$. We sometimes call $\alpha$ also an {\em
action\/} of $\calP$ on $V$.  }
\end{Appendix}

\begin{Appendix}{\rm
\label{odpalil_jsem_notebook}
The {\em suspension\/} $\ss A = 
\{\ss A(n)\}_{n \geq 0}$ of a $\Sigma$-module $A = \{A(n)\}_{n \geq 0}$
is defined by
\[
\ss A(n) := \hskip .2em \uparrow^{n-1} \hskip -.2em A(n) \otimes
\sgn_n,
\] 
where $\sgn_n$ denotes the signum representation of
$\Sigma_n$, see~\cite[Definition~II.3.15]{markl-shnider-stasheff:book}. 
If $\calP$ is an operad, then the collection 
$\ss \calP$ carries a canonical induced
operad structure and the operad $\ss \calP$ is called the {\em
operadic suspension\/} of $\calP$.
For any two operads $\calP$ and $\calQ$, 
\[
\ss (\calP \otimes \calQ) \cong \ss \calP \ot \calQ \cong \calP \ot
\ss \calQ. 
\] 
}
\end{Appendix}

\begin{Appendix}
\label{ma_prejit_fronta_v_sobotu}
{\rm 
An $(m,n)$-{\em algebra\/}
is~\cite[Example~9.4]{fox-markl:ContM97} a graded vector space $A$
together with two bilinear maps, $-\cup-:A \ot A \to A$ of degree $m$,
and $[-,-]:A\ot A\to A$ of degree $n$ ($m$ and $n$ are natural
numbers), such that, for any homogeneous $a,b,c\in A$,
\begin{itemize}
\item[(i)]
$a\cup b = (-1)^{|a|\cdot |b|+m}\cdot b\cup a$,
\item[(ii)]
$[a,b] = -(-1)^{|a|\cdot |b|+n}\cdot [b,a]$,
\item[(iii)]
$-\cup-$ is associative in the sense that
$$
a\cup (b\cup c) = (-1)^{m\cdot(|a|+1)}\cdot(a\cup b)\cup c,
$$
\item[(iv)]
$[-,-]$ satisfies the following form of the Jacobi identity:
$$
(-1)^{|a|\cdot(|c|+n)}\cdot [a,[b,c]]+(-1)^{|b|\cdot(|a|+n)}\cdot
[b,[c,a]]
+(-1)^{|c|\cdot(|b|+n)}\cdot [c,[a,b]] = 0,
$$
\item[(v)]
the operations $-\cup-$ and $[-,-]$ are compatible in the sense that
$$
(-1)^{m\cdot|a|}[a,b\cup c] =
[a,b]\cup c + (-1)^{(|b|\cdot |c|+m)}[a,c]\cup b.
$$
\end{itemize}
$(0,1)$-algebras were considered in~\cite{getzler-jones:preprint}
under the name {\em $2$-algebras\/} or {\em braid algebras\/}. The
corresponding operad $\Braid$ is isomorphic to the homology of the
little discs operad $\disc_2$, $\Braid \cong H_*(\disc_2)$.
Following~\cite[Section~10]{gerstenhaber-voronov:FAP95}, we call
$(1,0)$-algebras {\em Gerstenhaber algebras\/}, though the terminology
is not unique, compare for instance~\cite[Subsection~10.2]{ADT88}
where a Gerstenhaber algebra means a $(0,-1)$-algebra.
}
\end{Appendix}

\begin{Appendix}
\label{zase_do_Prahy}
{\rm
Let $M$ be a right module over a finite group $G$. We denote, as usual
\[
M^G := \{ m \in M;\ m g =g \mbox { for all } g \in G\}
\ \mbox { and }\
M_G := \frac M{(m-mg;\ m \in M,\ g \in G)}.
\]
Let ${\it Aver\/}: M \to M^G$ be the ``averaging'' map given by
\[
{\it Aver\/}(m):= \frac 1{|G|}\sum_{g \in G} mg.
\]
It is a standard fact that the composition $\pi \iota$ of the
projection $\pi : M \epi M_G$ with the inclusion $\iota : M^G
\hookrightarrow M$ is the identity and that ${\it Aver\/}$ is a right
inverse to $\iota$.
}
\end{Appendix}

\begin{Appendix}
\label{Ferda}
{\rm
Let us recall the operadic cochain complex
and introduce some useful notations.
As a graded vector space, the operadic cochain complex is defined
by~\cite[Definition~II.3.99]{markl-shnider-stasheff:book}:
\begin{equation}
\label{neco_tady_hrozne_kvika}
\CAP{n-1} = \left[ 
\Lin(\otexp {(\downarrow \hskip -.2em A)}n,\downarrow\hskip -.2em A) 
\ot \calP^!(n)
\right]^{\Sigma_n} \hskip -.9em, \ n \geq 1,
\end{equation}
where $\downarrow A$ denotes the desuspension of the graded vector
space $A$. 
It will be convenient to denote 
\[
\tildeCAP{n-1} := \Lin(\otexp {(\downarrow \hskip -.2em
  A)}n,\downarrow  \hskip -.2em A) \ot \calP^!(n),
\]
so that $\CAP {n-1} \cong \tildeCAP{n-1}^{\Sigma_n}$.
The averaging over the $\Sigma_n$-action defines an epimorphism
\[
{\it Aver\/} : \tildeCAP{n-1} \epi \CAP{n-1}
\]
of graded modules which is a left inverse to the inclusion
\[
\iota : \CAP{n-1} \hookrightarrow \tildeCAP{n-1}.
\] 

We will often use the following canonical isomorphisms of graded
$\Sigma_n$-modules:
\begin{eqnarray*}
\tildeCAP {n-1} &=&
\Lin(\otexp {(\downarrow \hskip -.2em A)}n,\downarrow\hskip -.2em A) 
\ot \calP^!(n) \cong
\hskip .2em \uparrow^{\ot {n-1}} \hskip -.2em 
(\Lin(\otexp An,A) \ot
\calP^!(n) \ot \sgn_n) 
\\
&\cong&
\ss(\Lin(\otexp An,A) \ot \calP^!(n))
\cong  \ss\End_A(n) \ot \calP^!(n) 
\\
&\cong& \End_A(n) \ot \ss\calP^!(n) \cong \End_{\desusp A}(n) \ot \calP^!(n).
\end{eqnarray*}
}
\end{Appendix}

\baselineskip 1.1em
\def\cprime{$'$}

\end{document}